\pdfoutput=1
\RequirePackage{ifpdf}
\ifpdf 
\documentclass[pdftex]{sigma}
\else
\documentclass{sigma}
\fi

\usepackage{amscd}
\usepackage{fancyhdr}
\usepackage{stmaryrd} 
\usepackage{bbm} 
\usepackage{mathtools}
\usepackage{tikz-cd}
\usepackage[arrow, matrix, curve]{xy}
\usepackage{extarrows}

\newcommand\circlearound[1]{
 \tikz[baseline]\node[draw,shape=circle,anchor=base] {#1} ;}

\newcommand{\R}{\mathcal{R}}

\DeclareMathOperator{\supp}{supp}

\numberwithin{equation}{section}

\newtheorem{Theorem}{Theorem}[section]
\newtheorem*{Theorem*}{Theorem}
\newtheorem{Corollary}[Theorem]{Corollary}
\newtheorem{Lemma}[Theorem]{Lemma}
\newtheorem{Proposition}[Theorem]{Proposition}

\theoremstyle{definition}
\newtheorem{Definition}[Theorem]{Definition}

\newtheorem{Example}[Theorem]{Example}
\newtheorem{Remark}[Theorem]{Remark}
\newtheorem*{Remark*}{Remark}

\begin{document}
\allowdisplaybreaks

\newcommand{\arXivNumber}{2502.02353}

\renewcommand{\PaperNumber}{020}

\FirstPageHeading

\ShortArticleName{Motives of Nullcones of Quiver Representations}

\ArticleName{Motives of Nullcones of Quiver Representations}

\Author{Lydia G\"OSMANN and Markus REINEKE}
\AuthorNameForHeading{L. G\"osmann and M. Reineke}
\Address{Faculty of Mathematics, Ruhr University Bochum, D-44801 Bochum, Germany}
\Email{\mail{lydia.goesmann@rub.de}, \mail{markus.reineke@rub.de}}

\ArticleDates{Received June 13, 2025, in final form February 10, 2026; Published online March 02, 2026}

\Abstract{We derive two formulas for motives of nullcones of quiver representations, one being explicit, the other being of wall-crossing type.}


\Classification{14C15; 14L24; 16G20}

\section{Introduction}

For a finite-dimensional complex representation of a complex reductive group, its nullcone, already introduced by Hilbert \cite{H}, plays an important role in the study of the invariant ring, the orbit structure, and the quotient map for the group action \cite{KW}.

A particularly interesting and accessible class of representations are the representation spaces of quivers with their base change action of a product of general linear groups. In this case, the orbits correspond to isomorphism classes of quiver representations, and the quotient realises a~moduli space for isomorphism classes of semisimple representations \cite{LP}.

It is thus natural, both from the invariant theoretic and from the quiver-theoretic point of view, to study the nullcones of representation spaces of quivers, which parametrize nilpotent representations. A qualitative study of the nullcone was accomplished in \cite{LB}, adapting the Hesselink stratification of the nullcone, and describing it in terms of loci of semistable representations for so-called level quivers.

In the present work, we perform a systematic quantitative study of these nullcones, by describing their motives, that is, their classes in the Grothendieck ring of complex varieties. A first such formula (more precisely, a formula for counts of rational points over finite fields) already appears in \cite[Lemma 2.6]{BSV}.

Our first result gives a simple explicit formula for the motives, as a recursion formula (Corollary~\ref{cor31}), a generating function identity (Corollary~\ref{gfi}), and as an explicit sum formula (Corollary~\ref{resolved}). These formulas are derived from a geometric construction essentially using the stratification of nilpotent quiver representations by their socle type, inspired by methods of \cite{RCounting} (Section~\ref{s2}).

This result can be seen as a generalization of a prototypical motivic description of nullcones, namely the classical result \cite{FH} that there are $q^{d(d-1)}$ nilpotent $d\times d$ matrices over a finite field with $q$ elements.

Our second result (Sections \ref{s4} and \ref{s5}) interprets the Hesselink stratification of~\cite{LB} as an identity of motives, first yielding a positive summation formula (Theorem~\ref{theo mot}) for the motive of the nullcone, and culminating in a product formula relating the generating function of motives of nullcones to generating functions of semistable representations of level quivers (Theorem~\ref{theoremwcf}), reminiscent of wall-crossing formulas such as \cite[Theorem~2.4]{MozIntro}.

It is reasonable to expect that our formulas also hold in positive characteristic, and thus yield formulas for numbers of rational points over finite fields. However, the methods of \cite{B,Hesse} which we use in an essential way require a base field of characteristic zero.

\section{An explicit recursive formula for the motive of the nullcone}\label{s2}

\subsection{Quiver representations and their nullcone}

Let $Q$ be a finite quiver with set of vertices $Q_0$ and arrows written $\alpha\colon i\rightarrow j$. Define the Euler form $\langle\_,\_\rangle$ of $Q$ by
\[\langle{\bf d},{\bf e}\rangle=\sum_{i\in Q_0}d_ie_i-\sum_{\alpha\colon i\rightarrow j}d_ie_j\]
for ${\bf d}=(d_i)_i$, ${\bf e}=(e_i)_i\in\mathbb{Z}Q_0$. We will consider complex representations of $Q$ of some dimension vector $\mathbf{d}\in\mathbb{N}Q_0$. In particular, we denote by $S_i$ the one-dimensional representation given by~$(S_i)_i=\mathbb{C}$, and all arrows being represented by zero maps.

Fixing vector spaces $V_i$ of dimension $d_i$ for $i\in Q_0$, respectively, let \[R_\mathbf{d}=R_\mathbf{d}(Q)=\bigoplus_{\alpha\colon i\rightarrow j}{\rm Hom}_\mathbb{C}(V_i,V_j)\] be the space of representations of $Q$ on the $V_i$, acted on by the base change group \[G_\mathbf{d}=\prod_{i\in Q_0}{\rm GL}(V_i)\] by $Q$-graded conjugation
$(g_i)_i\cdot(V_\alpha)_\alpha=\bigl(g_jV_\alpha g_i^{-1}\bigr)_{\alpha\colon i\rightarrow j}$.

Let $N_\mathbf{d}$ be the nullcone for this action. By its definition and the Hilbert criterion, it can be characterized by the following equivalent conditions.

\begin{Theorem}[{\cite[Section 1.4]{LB}}] The following are equivalent for a representation $V\in R_\mathbf{d}$:
\begin{enumerate}\itemsep=0pt
\item[$(1)$] $V$ belongs to $N_\mathbf{d}$,
\item[$(2)$] $0$ belongs to the closure of the orbit $G_{\bf d}\cdot V$,
\item[$(3)$] all $G_{\bf d}$-invariants on $R_{\bf d}$ vanish at $V$, that is $($by {\rm\cite{LP})}, ${\rm tr}(V_{\alpha_1}\circ\dots\circ V_{\alpha_s})=0$ for all oriented cycles $\alpha_1\cdots\alpha_s$ in $Q$,
\item[$(4)$] all operators $V_{\alpha_1}\circ\dots\circ V_{\alpha_s}$ for oriented cycles $\alpha_1\cdots\alpha_s$ in $Q$ are nilpotent,
\item[$(5)$] $V$ admits a filtration by subrepresentations with subquotients $S_i$ for $i\in Q_0$ {\rm (}the one-dimensional simple representations of $Q)$.
\end{enumerate}
\end{Theorem}

\subsection{Motives}

We denote by $K_0({\rm Var}_\mathbb{C})$ the Grothendieck ring of complex varieties; the class of a variety $X$ therein is denoted by $[X]$. We denote by $\mathbb{L}=\big[\mathbb{A}^1\big]$ the Lefschetz motive. We will work with the localized base ring of motives \[\mathbb{M}=K_0({\rm Var}_\mathbb{C})\big[\mathbb{L}^{-1},\bigl(1-\mathbb{L}^i\bigr)^{-1} : i\geq 1\big]\]
(the role of the localization is a purely formal one, to have the flexibility of dividing by motives of general linear groups). We define the Pochhammer symbol $(q)_n=(1-q)\cdots(1-q^n)$. Using this notation, we have the following formulas for the motives of general linear groups and Grassmannians
\[[{\rm GL}_n(\mathbb{C})]=(-1)^n\mathbb{L}^{n(n-1)/2}(\mathbb{L})_n,\qquad [{\rm Gr}_k^n]=\left[{n\atop k}\right]_\mathbb{L}:=\frac{(\mathbb{L})_n}{(\mathbb{L})_k\cdot(\mathbb{L})_{n-k}}.\]

We also note that, for varieties $X$ and $Y$ admitting a bijective morphism $f\colon X\rightarrow Y$, we have~${[X]=[Y]}$ in $K_0({\rm Var}_\mathbb{C})$ \cite[Lemma~2.9]{B}.

\subsection{The recursive formula}

Our aim is to describe various explicit formulas for the motive $[N_{\bf d}]\in \mathbb{M}$. We introduce some shorthand notation concerning dimension vectors; namely, we define
\begin{gather*}
\dim {\bf d}=\sum_id_i, (q)_{\bf d}=\prod_{i\in Q_0}(q)_{d_i}, \left[{{\bf d}\atop{\bf e}}\right]_q=\prod_{i}\left[{d_i\atop e_i}\right]_q,\\
 \left[{{\bf d}\atop{{\bf d}^1,\dots,{\bf d}^s}}\right]_q=\prod_i\frac{(q)_{d_i}}{\prod_{k=1}^s(q)_{d^k_i}}\end{gather*}
for ${\bf d}={\bf d}^1+\dots+{\bf d}^s\in\mathbb{N}Q_0$.

\begin{Theorem}\label{thm1} We have the following recursion determining $[N_{\bf d}]$ for all ${\bf d}\not=0$:
\[\frac{[N_{\bf d}]}{[G_{\bf d}]}=-\sum_{{\bf e}\lneqq{\bf d}}\mathbb{L}^{-\langle{\bf e},{\bf d}-{\bf e}\rangle}\cdot(\mathbb{L})_{{\bf d}-{\bf e}}^{-1}\cdot\frac{[N_{\bf e}]}{[G_{\bf e}]}.\]
\end{Theorem}

\begin{proof} Assume ${\bf d}\not=0$. For a vertex $i\in Q_0$ and a representation $V\in N_{\bf d}$ given by a~tuple~$(V_\alpha)_\alpha$, we consider the map
\[\Phi_i(V)=\bigoplus_{\alpha\colon i\rightarrow j}V_\alpha\colon\ V_i\rightarrow\bigoplus_{\alpha\colon i\rightarrow j}V_j.\]

We define a stratification of $N_{\bf d}$ by locally closed subsets $S_{\bf f}$ for ${\bf f}\leq{\bf d}$ by
\[S_{\bf f}=\{(V_\alpha)_\alpha \mid \dim{\rm Ker}(\Phi_i(V))=f_i\mbox{ for all }i\in Q_0\}.\]

Since a nilpotent representation is an iterated extension of the one-dimensional simple representations $S_i$ of $Q$, it admits some $S_i$ as a subrepresentation, and thus belongs to some $S_{\bf f}$ for~${{\bf f}\not={ 0}}$. In other words, $V$ belongs to $S_{\bf f}$ if and only if the socle of $V$ is of dimension vector ${\bf f}$, which is nonzero since $V$ is so. In particular, we have $[S_{ 0}]=0$.

For ${\bf f}\leq{\bf d}$, we consider the product of Grassmannians
\[{\rm Gr}_{\bf f}({\bf d})=\prod_{i\in Q_0}{\rm Gr}_{f_i}(V_i),\]
and, more generally, for ${\bf f}\leq{\bf g}\leq{\bf d}$, the corresponding product \smash{${\rm Fl}_{{\bf f},{\bf g}}({\bf d})$} of two-step flag varieties.
Inside ${\rm Gr}_{\bf f}({\bf d})\times N_{\bf d}$, we consider the closed subvariety
\[X_{\bf f}=\{((U_i)_i,(V_\alpha)_\alpha) \mid U_i\subset{\rm Ker}(\Phi_i(V))\mbox{ for all }i\in Q_0\}.\]
The first projection $p_1\colon X_{\bf f}\rightarrow {\rm Gr}_{\bf f}({\bf d})$ being equivariant for the obvious $G_{\bf d}$-actions, we have $X_{\bf f}\simeq G_{\bf d}\times^{P_{\bf f}}Z_{\bf f}$,
where, for a fixed choice of subspaces $\bigl(U_i^0\bigr)_i\in {\rm Gr}_{\bf f}({\bf d})$, the group $P_{\bf f}\subset G_{\bf d}$ is the maximal parabolic fixing the $U_i^0$, and $Z_{\bf f}$ is the closed subset of $N_{\bf d}$ of representations for which every~$V_\alpha$ for $\alpha\colon i\rightarrow j$ annihilates $U_i^0$. Choosing bases of the $V_i$ compatible with the $U_i^0$, in such representations the arrows of $Q$ are represented by block matrices
\[
\left[\begin{matrix}0&\zeta_\alpha\\ 0&W_\alpha\end{matrix}\right],
\]
where $(W_\alpha)$ defines a representation in $N_{{\bf d}-{\bf f}}$. It follows that $Z_{\bf f}$ is a trivial vector bundle of rank \smash{$\sum_{\alpha\colon i\rightarrow j}(d_i-f_i)f_j$} over $N_{{\bf d}-{\bf f}}$. Consequently, we find for the motive of $X_{\bf f}$
\[[X_{\bf f}]=\frac{[G_{\bf d}]}{[P_{\bf f}]}\cdot[Z_{\bf f}]=\frac{[G_{\bf d}]}{[G_{{\bf d}-{\bf f}}]\cdot[G_{\bf f}]}\cdot\mathbb{L}^{-\langle{\bf d}-{\bf f},{\bf f}\rangle}\cdot[N_{{\bf d}-{\bf f}}].\]

On the other hand, we consider the second projection $p_2\colon X_{\bf f}\rightarrow N_{\bf d}$. By definition of the strata $S_{\bf f}$, the image of $p_2$ is the union of the $S_{\bf g}$ for ${\bf f}\leq{\bf g}\leq{\bf d}$. We denote the inverse image of such a stratum $S_{\bf g}$ under $p_2$ by $X_{{\bf f},{\bf g}}$. By definition, the projection
$p_2\colon X_{{\bf f},{\bf g}}\rightarrow S_{\bf g}$
has all its fibres isomorphic to ${\rm Gr}_{\bf f}({\bf g})$. We claim that we have
$[X_{{\bf f},{\bf g}}]=[S_{\bf g}]\cdot[{\rm Gr}_{\bf f}({\bf g})]$
in $\mathbb{M}$. For this, we observe that we have a diagram of Cartesian squares
\[
\begin{array}{ccccc}{\rm Fl}_{{\bf f},{\bf g}}({\bf d})&\leftarrow&\widetilde{X}&\stackrel{\widetilde{p}}{\rightarrow}&X_{{\bf f},{\bf g}}\\\downarrow&&\downarrow&&\downarrow\\
{\rm Gr}_{\bf g}({\bf d})&\leftarrow&X_{{\bf g},{\bf g}}&\stackrel{p}{\rightarrow}&S_{\bf g}.\end{array}
\]

Namely, from both squares the variety $\widetilde{X}$ is calculated as
\[\widetilde{X}=\{((U_i)_i,(U_i')_i,(V_\alpha)_\alpha) \mid U_i\subset U_i'={\rm Ker}(\Phi_i(V))\mbox{ for all }i\in Q_0\}.\]

As a special case of the previous observation, we see that the map $p\colon X_{{\bf g},{\bf g}}\rightarrow S_{\bf g}$ is a bijection (but we do not know a priori whether it is an isomorphism), thus $\widetilde{p}$ is a bijection, too. The map~${{\rm Fl}_{{\bf f},{\bf g}}({\bf d})\rightarrow{\rm Gr}_{\bf g}({\bf d})}$ is Zariski-locally trivial due to $G_{\bf d}$-equivariance, thus the same holds for the map $\widetilde{X}\rightarrow X_{{\bf g},{\bf g}}$. It follows that we can calculate
\[[X_{{\bf f},{\bf g}}]=\big[\widetilde{X}\big]=[X_{{\bf g},{\bf g}}]\cdot[{\rm Gr}_{\bf f}({\bf g})]=[S_{\bf g}]\cdot[{\rm Gr}_{\bf f}({\bf g})],\]
as claimed.\\[1ex]
Consequently, we have
\[[X_{\bf f}]=\sum_{{\bf f}\leq{\bf g}\leq{\bf d}}[S_{\bf g}]\cdot[{\rm Gr}_{\bf f}({\bf g})].\]

We note the following identity of Gaussian binomial coefficients for all $m\geq 0$ (see, for example, \cite[equation~(3.14)]{RCounting})
\[\sum_{a=0}^m(-1)^a\mathbb{L}^{a(a-1)/2}\left[{m\atop a}\right]_\mathbb{L}=\delta_{m,0}.\]

This implies the following identity for all ${\bf d}\not=0$:
\begin{align*}
\sum_{{\bf f}\leq{\bf d}}(-1)^{\dim{\bf f}}\mathbb{L}^{\sum_i f_i(f_i-1)/2}[X_{\bf f}]&=\sum_{{\bf f}\leq{\bf g}\leq{\bf d}}(-1)^{\dim{\bf f}}\mathbb{L}^{\sum_i f_i(f_i-1)/2}[{\rm Gr}_{\bf f}({\bf g})]\cdot[S_{\bf g}]
\\
&=\sum_{{\bf g}\leq{\bf d}}\underbrace{\bigg(\sum_{{\bf f}\leq{\bf g}}(-1)^{\dim {\bf f}}\mathbb{L}^{\sum_if_i(f_i-1)/2}[{\rm Gr}_{\bf f}({\bf g})]\bigg)}_{=\delta_{{\bf g},{0}}}\cdot[S_{\bf g}]\\
&=[S_{ 0}]=0.
\end{align*}
Substituting the above formula for $[X_{\bf f}]$, we find
\[0=\sum_{{\bf f}\leq{\bf d}}(-1)^{\dim{\bf f}}\mathbb{L}^{\sum_i f_i(f_i-1)/2}\frac{[G_{\bf d}]}{[G_{{\bf d}-{\bf f}}]\cdot[G_{\bf f}]}\cdot\mathbb{L}^{-\langle{\bf d}-{\bf f},{\bf f}\rangle}\cdot[N_{{\bf d}-{\bf f}}].\]
Dividing by $[G_{\bf d}]$ and isolating the term corresponding to ${\bf f}={ 0}$, we find
\[\frac{[N_{\bf d}]}{[G_{\bf d}]}=-\sum_{{0}\not={\bf f}\leq{\bf d}}\underbrace{\frac{(-1)^{\sum_if_i}\mathbb{L}^{\sum_if_i(f_i-1)/2}}{[G_{\bf f}]}}_{=(\mathbb{L})_{{\bf f}}^{-1}}\cdot\mathbb{L}^{-\langle{\bf d}-{\bf f},{\bf f}\rangle}\cdot\frac{[N_{{\bf d}-{\bf f}}]}{[G_{{\bf d}-{\bf f}}]}.\]
Finally, substituting ${\bf f}={\bf d}-{\bf e}$, the claimed formula follows, finishing the proof.\end{proof}

\section{Applications and examples}\label{s3}

\subsection{Applications}

We will now reformulate the recursive formula of Theorem~\ref{thm1} in several ways.

\begin{Corollary}\label{cor31} We have the following equivalent form of the recursion, which holds as an identity in $K_0({\rm Var}_\mathbb{C})$:
\[[N_{\bf d}]=-\sum_{{\bf e}\lneqq{\bf d}}(-1)^{\dim{\bf d}-\dim{\bf e}}\cdot\mathbb{L}^{\sum_i\big({d_i\choose 2}-{e_i\choose 2}\big)-\langle{\bf e},{\bf d}-{\bf e}\rangle}\cdot\left[{{\bf d}\atop{\bf e}}\right]_\mathbb{L}\cdot{[N_{\bf e}]}.\]
In particular, there exist polynomials $n_{\bf d}(x)\in\mathbb{Z}[x]$ such that $[N_{\bf d}]=n_{\bf d}(\mathbb{L})$.
\end{Corollary}

\begin{proof} The recursion follows from Theorem~\ref{thm1} just by making the motive of $G_{\bf e}$ explicit. Defining polynomials $n_{\bf d}$ by the same recursion with $\mathbb{L}$ replaced by $x$, and starting with $n_0=1$, we obviously have $n_{\bf d}\in\mathbb{Z}\big[x^{\pm 1}\big]$, and $n_{\bf d}(\mathbb{L})=[N_{\bf d}]$. By the definition of the Euler form, the $\mathbb{L}$-exponent in the recursion is always nonnegative, thus we have $n_{\bf d}(x)\in\mathbb{Z}[x]$.\end{proof}

We define the ring $\mathbb{M}_\mathbb{L}\llbracket Q_0\rrbracket$ as the twisted formal power series ring with topological $\mathbb{M}$-basis elements $t^{\bf d}$ for ${\bf d}\in \mathbb{N}Q_0$ and multiplication
\smash{$t^{\bf d}\cdot t^{\bf e}=\mathbb{L}^{-\langle{\bf d},{\bf e}\rangle}\cdot t^{{\bf d}+{\bf e}}$}.
Such twists are frequently used for identities of generating series of motives, see, for example, \cite[Section~2.1]{MozIntro} and \cite[Definition~3.2]{RCounting}.

\begin{Corollary}\label{gfi} In $\mathbb{M}_\mathbb{L}\llbracket Q_0\rrbracket$, we have the identity
\[\sum_{\bf d}\frac{[N_{\bf d}]}{[G_{\bf d}]}t^{\bf d}\cdot\sum_{\bf d}\frac{t^{\bf d}}{(\mathbb{L})_{\bf d}}=1.\]
\end{Corollary}

\begin{proof} Expanding the left-hand side of the claimed equation, its $t^{\bf d}$-coefficient equals
\[\sum_{{\bf e}\leq{\bf d}}\mathbb{L}^{-\langle{\bf e},{\bf d}-{\bf e}\rangle}\cdot(\mathbb{L})_{{\bf d}-{\bf e}}^{-1}\cdot\frac{[N_{\bf e}]}{[G_{\bf e}]},\]
which equals $0$ for ${\bf d}\not=0$ by Theorem~\ref{thm1}.\end{proof}

Using this identity of generating functions, we can resolve the recursion of Theorem~\ref{thm1}.

\begin{Corollary}\label{resolved} We have
\[[N_{\bf d}]=(-1)^{\dim {\bf d}}\cdot\mathbb{L}^{\sum_i{d_i\choose 2}}\cdot\sum_{{\bf d}^*}(-1)^s\cdot\mathbb{L}^{-\sum_{k<l}\langle{\bf d}^k,{\bf d}^l\rangle}\cdot\left[{{\bf d}\atop{{\bf d}^1,\dots,{\bf d}^s}}\right]_\mathbb{L},\]
where the sum ranges over all ordered decompositions ${\bf d}={\bf d}^1+\dots+{\bf d}^s$ of ${\bf d}$ into non-zero parts~${\bf d}^k$.
\end{Corollary}

\begin{proof} By Corollary~\ref{gfi}, we can calculate
\begin{align*}
\sum_{\bf d}\frac{[N_{\bf d}]}{[G_{\bf d}]}t^{\bf d}&=\biggl(1+\sum_{\bf d\not={0}}\frac{t^{\bf d}}{(\mathbb{L})_{\bf d}}\biggr)^{-1}=\sum_{s\geq 0}(-1)^s\biggl(\sum_{{\bf d}\not=0}\frac{t^{\bf d}}{(\mathbb{L})_{\bf d}}\biggr)^s\\
&=\sum_{s\geq 0}(-1)^s\cdot
\sum_{{\bf d}^1,\dots,{\bf d}^s\not=0}\mathbb{L}^{-\sum_{k<l}\langle{\bf d}^k,{\bf d}^l\rangle}
\frac{t^{\bf d}}{(\mathbb{L})_{{\bf d}^1}\cdots(\mathbb{L})_{{\bf d}^s}}.
\end{align*}
Comparing coefficients of $t^{\bf d}$ on both sides and rewriting $[G_{\bf d}]$ using the Pochhammer symbol, we arrive at the claimed formula.\end{proof}

We call a quiver symmetric if its Euler form is symmetric. Equivalently, denoting by $r_{ij}$ the number of arrows from $i$ to $j$ in $Q$ for $i,j\in Q_0$, symmetry is given by $r_{ij}=r_{ji}$ for all~${i,j\in Q_0}$. In this case, it is immediately verified that
\smash{$x^{\bf d}=\bigl(-\mathbb{L}^{1/2}\bigr)^{-\langle{\bf d},{\bf d}\rangle}t^{\bf d}$} fulfills $x^{\bf d}\cdot x^{\bf e}=x^{{\bf d}+{\bf e}}$; consequently, we do not have to work with twisted multiplications of formal series in the symmetric case. We extend the base ring of motives by adjoinging a square root of $\mathbb{L}$. This allows us to define motivic Donaldson--Thomas invariants (see \cite{RDT} for an exposition) by the identity
\[\sum_{\mathbf d}\bigl(-\mathbb{L}^{1/2}\bigr)^{\langle\mathbf{d},\mathbf{d}\rangle}\frac{[R_\mathbf{d}]}{[G_\mathbf{d}]}x^\mathbf{d} =\operatorname{Exp}\biggl(\frac{1}{\mathbb{L}^{-1/2}-\mathbb{L}^{1/2}}\sum_{\mathbf{d}\not=0}{\rm DT}^Q_\mathbf{d}(\mathbb{L}) x^\mathbf{d}\biggr)\]
in $\mathbb{M}\big[\mathbb{L}^{1/2}\big]\llbracket Q_0\rrbracket$. Here $\operatorname{Exp}$ is defined by
\[\operatorname{Exp}\biggl(\sum_{{\bf d}\not=0}\sum_kc_{{\bf d},k}q^{k/2}x^{\bf d}\biggr)=\prod_{{\bf d}\not=0}\prod_k\bigl(1-q^{k/2}x^{\bf d}\bigr)^{-c_{{\bf d},k}}.\]
It is known by \cite{E} that the Donaldson--Thomas invariant is a polynomial in $-\mathbb{L}^{1/2}$ with nonnegative coefficients, that is, ${\rm DT}^Q_\mathbf{d}(\mathbb{L})\in \mathbb{Z}\big[\mathbb{L}^{\pm 1/2}\big]$.

\begin{Corollary} Assume that $Q$ is symmetric. Then \[\sum_\mathbf{d}\bigl(-\mathbb{L}^{1/2}\bigr)^{\langle\mathbf{d},\mathbf{d}\rangle} \frac{[N_\mathbf{d}]}{[G_\mathbf{d}]}x^\mathbf{d}=\operatorname{Exp}\biggl(\frac{1}{\mathbb{L}^{-1/2}-\mathbb{L}^{1/2}}\sum_{\mathbf{d}\not=0}{\rm DT}^Q_{\bf d}\bigl(\mathbb{L}^{-1}\bigr)x^\mathbf{d}\biggr).\]
\end{Corollary}

\begin{proof} The left-hand side of the identity defining the DT invariants equals
\[\sum_{\bf d}\bigl(-\mathbb{L}^{1/2}\bigr)^{-\langle{\bf d},{\bf d}\rangle}\frac{x^{\bf d}}{\bigl(\mathbb{L}^{-1}\bigr)_{\bf d}}.\]
Polynomiality of the Donaldson--Thomas invariants allows us to formally interchange $\mathbb{L}$ and $\mathbb{L}^{-1}$, and thus we have
\[\sum_{\bf d}\bigl(-\mathbb{L}^{1/2}\bigr)^{\langle{\bf d},{\bf d}\rangle}\frac{x^{\bf d}}{(\mathbb{L})_{\bf d}}=\operatorname{Exp}\biggl(\frac{1}{\mathbb{L}^{1/2}-\mathbb{L}^{-1/2}}\sum_{{\bf d}\not=0}{\rm DT}^Q_{\bf d}\bigl(\mathbb{L}^{-1}\bigr)x^{\bf d}\biggr).\]
On the other hand, using the definition of $x^{\bf d}$ and Corollary~\ref{gfi}, we have
\[\sum_\mathbf{d}\bigl(-\mathbb{L}^{1/2}\bigr)^{\langle\mathbf{d},\mathbf{d}\rangle} \frac{[N_\mathbf{d}]}{[G_\mathbf{d}]}x^\mathbf{d}=\biggl(\sum_{\bf d}\bigl(-\mathbb{L}^{1/2}\bigr)^{\langle{\bf d},{\bf d}\rangle}\frac{x^{\bf d}}{(\mathbb{L})_{\bf d}}\biggr)^{-1}.\]
Combining these identities, the claim follows.\end{proof}

\begin{Corollary} If $Q$ is symmetric with $r_{ij}$ arrows from $i$ to $j$ for $i,j\in Q_0$, the polynomial~$n_{\bf d}(x)$ has leading term
\[\frac{(\dim {\bf d})!}{\prod_{i\in Q_0}d_i!}\cdot x^{\sum_{i}(r_{ii}+1){d_i\choose 2}+\sum_{i<j}r_{ij}d_id_j}.\]
\end{Corollary}

\begin{proof} By Corollary~\ref{resolved}, we have
\[n_{\bf d}(x)=(-1)^{\dim {\bf d}}\cdot{x}^{\sum_i{d_i\choose 2}}\cdot\sum_{{\bf d}^*}(-1)^s\cdot{x}^{-\sum_{k<l}\langle{\bf d}^k,{\bf d}^l\rangle}\cdot\left[{{\bf d}\atop{{\bf d}^1,\dots,{\bf d}^s}}\right]_{x}.\]
The degree of the summand corresponding to a decomposition ${\bf d}^*$ is, by symmetry of $Q$, given~by
\[-\sum_{k<l}\big\langle{\bf d}^k,{\bf d}^l\big\rangle+\sum_{k<l}\sum_id^k_id^l_i=\sum_{k<l}\sum_{i,j}r_{ij}d^k_id^l_j=\frac{1}{2} \biggl(\sum_{i,j}r_{ij}d_id_j-\sum_{i,j}\sum_kr_{ij}d^k_id^k_j\biggr).\]

Refining an ordered partition ${\bf d}^*$ by replacing a part ${\bf d}^k$ by two parts ${\bf d}'$, ${\bf d}''$ (in this order) thus increases the degree of the corresponding summand by $\sum_{i,j}r_{ij}d'_id''_j$. The degree can thus be maximal only for maximally refined decompositions, which are of the form ${\bf i}_1,\dots,{\bf i}_s$, where each ${\bf i}$ for $i\in Q_0$ appears $d_i$-times in the sequence; in particular, $s=\dim {\bf d}$. We define \[K_i=\{k=1,\dots,s \mid {\bf i}_k=i\}\] for $i\in Q_0$. The summand corresponding to the sequence ${\bf i}_1,\dots,{\bf i}_s$ is then of degree
\begin{align*}
\frac{1}{2}\biggl(\sum_{i,j}r_{ij}d_id_j-\sum_{i,j}\sum_kr_{ij}\delta_{k\in K_i}\delta_{k\in K_j}\biggr)
&=\frac{1}{2}\biggl(\sum_{i,j}r_{ij}d_id_j-\sum_i\sum_kr_{ii}\biggr)\\
&=\sum_ir_{ii}{d_i\choose 2}+\sum_{i<j}r_{ij}d_id_j.\end{align*}
All of the corresponding summands have the same sign $(-1)^{\dim {\bf d}}$, thus there are no cancellations in the summation of these top degree terms. The number of these summands equals the number of decompositions of $\{1,\dots,\dim{\bf d}\}$ into parts $K_i$ for $i\in Q_0$, which is the quotient of factorials of the claimed formula. Adding the contribution $\sum_i{d_i\choose 2}$ of the prefactor in the above formula for $n_{\bf d}(x)$ to the degree, we arrive at the claimed formula for the leading term.\end{proof}

\begin{Remark}
A similar proof appears in \cite[Theorem~4.1]{RSmall}. Refining the argument there shows that, in case $Q$ is symmetric, the variety $N_{\bf d}$ is equidimensional of dimension
\[\sum _{i}(r_{ii}+1){d_i\choose 2}+\sum_{i<j}r_{ij}d_id_j,\]
with ${(\dim{\bf d})!}/{\prod_id_i!}$
irreducible components.
\end{Remark}

\subsection{Examples}\label{ex part 1}

We finish this section with several examples illustrating the previous formulas.
\begin{enumerate}\itemsep=0pt
\item[(1)] If $Q$ is acyclic, every representation is nilpotent, thus
\smash{$[N_{\bf d}]=[R_{\bf d}]=\mathbb{L}^{\sum_{\alpha\colon i\rightarrow j}d_id_j}$}
trivially, but this is not obvious from any of the previous formulas. We will use two variants of the $q$-binomial theorem to derive this fact from Corollary~\ref{gfi}. Namely, we have
\[\sum_{d\geq 0}\frac{t^d}{(q)_d}=\prod_{k\geq 0}\bigl(1-q^kt\bigr)^{-1}, \sum_{d\geq 0}\frac{q^{d\choose 2}t^d}{(q)_d}=\prod_{k\geq 0}\bigl(1+q^kt\bigr),\]
and thus
\[\left(\sum_{d\geq 0}\frac{q^{d\choose 2}t^d}{(q)_d}\right)^{-1}=\prod_{k\geq 0}\bigl(1+q^kt\bigr)^{-1}=\sum_{d\geq 0}\frac{(-t)^d}{(q)_d}.\]
Since $Q$ is acyclic, we can enumerate its vertices $Q_0=\{i_1,\dots,i_n\}$ such that the existence of an arrow $i_k\rightarrow i_l$ implies $k<l$. By the definition of the twisted multiplication in~$\mathbb{M}_\mathbb{L}\llbracket Q_0\rrbracket$, we then have
\[t_{i_1}^{d_{i_1}}\cdots t_{i_n}^{d_{i_n}}=\mathbb{L}^{-\langle{\bf d},{\bf d}\rangle+\sum_kd_{i_k}(d_{i_k}+1)/2}t^{\bf d},\qquad t_{i_n}^{d_{i_n}}\cdots t_{i_1}^{d_{i_1}}=\mathbb{L}^{-\sum_kd_{i_k}(d_{i_k}-1)/2}t^{\bf d}.\]
We can thus rewrite
\[\sum_{\bf d}\frac{t^{\bf d}}{(\mathbb{L})_{\bf d}}=\prod_{k=1}^n\sum_{d\geq 0}\frac{\mathbb{L}^{d(d-1)/2}t_{i_{n+1-k}}^d}{(\mathbb{L})_d}\]
and then, by the above identity,
\begin{align*}
\left(\sum_{\bf d}\frac{t^{\bf d}}{(\mathbb{L})_{\bf d}}\right)^{-1}&=\prod_{k=1}^n\sum_{d\geq 0}\frac{(-t_{i_k})^d}{(\mathbb{L})_d}
=\sum_{\bf d}(-1)^{\dim{\bf d}}\frac{\mathbb{L}^{-\langle{\bf d},{\bf d}\rangle+\sum_id_i(d_i+1)/2}}{(\mathbb{L})_{\bf d}}t^{\bf d}\\
&=\sum_{\bf d}\frac{[R_{\bf d}]}{[G_{\bf d}]}t^{\bf d}.\end{align*}
Now $[N_{\bf d}]=[R_{\bf d}]$ follows from Corollary~\ref{gfi}.
\item[(2)] Similarly, if $Q$ is the quiver with a single vertex and one loop, Corollary~\ref{gfi}, together with the above two versions of the $q$-binomial theorem, easily yields
$[N_d]=\mathbb{L}^{d(d-1)}$
for $N_d$ the variety of nilpotent $d\times d$ matrices. This provides a motivic analogue of the classical result~\cite{FH}.
\item[(3)] If $Q$ is arbitrary and $d_i=1$ for all $i\in Q_0$, Corollary~\ref{resolved} reads
\[[N_{\bf 1}]=(-1)^{Q_0}\cdot\sum_{I_*}(-1)^s\cdot\mathbb{L}^{\sum_{k<l}|I_k\rightarrow I_l|},\]
where the sum ranges over all tuples $I_*=(I_1,\dots,I_s)$ of subsets of $Q_0$ inducing a disjoint union $Q_0=\bigcup_{k=1}^sI_k$, and $|I_k\rightarrow I_l|$ denotes the number of arrows from vertices in $I_k$ to vertices in $I_l$. On the other hand, denote by $Q_1$ the set of arrows in $Q$, and define the support of a representation $V\in R_{\bf 1}$ as the subset of $Q_1$ of all arrows $\alpha$ such that $V_\alpha\not=0$. Then $V$ is nilpotent if and only if its support is an acyclic subset $F\subset Q_1$, that is, such that the corresponding subquiver is acyclic. For any such acyclic subset, the locally closed subset of $R_{\bf 1}$ of representations having support $F$ is isomorphic to a torus of rank $|F|$, thus has motive $(\mathbb{L}-1)^{|F|}$ (as an aside, we note that every irreducible component of $N_{\bf 1}$ is thus toric). We then arrive at the formula
\smash{$[N_{\bf 1}]=\sum_{{F\subset Q_1}\atop\mbox{\tiny acyclic}}(\mathbb{L}-1)^{|F|}$}.
\item[(4)] For the quiver $Q$ with a single vertex and $m$ loops, the variety $N_{\bf d}$ consists of tuples of~$d\times d$ matrices $(A_1,\dots,A_m)$ such that every polynomial in these matrices is nilpotent, or, equivalently, which can be simultaneously conjugated to strictly upper triangular matrices. We list $[N_d]$ for small $d$ using Corollary~\ref{resolved}
\begin{gather*}
[N_1]=1, \qquad [N_2]=\mathbb{L}^{m+1}+\mathbb{L}^m-\mathbb{L},\\
[N_3]=\mathbb{L}^{3m+3}+2\mathbb{L}^{3m+2}+2\mathbb{L}^{3m+1}+\mathbb{L}^{3m}
-2\mathbb{L}^{2m+3}-2\mathbb{L}^{2m+2}-2\mathbb{L}^{2m+1}+\mathbb{L}^3,\\
[N_4]=\mathbb{L}^{6m+6}+3\mathbb{L}^{6m+5}+5\mathbb{L}^{6m+4}+6\mathbb{L}^{6m+3}+5\mathbb{L}^{6m+2}+3\mathbb{L}^{6m+1}+\mathbb{L}^{6m}\\
\phantom{[N_4]=}{}-3\mathbb{L}^{5m+6}-6\mathbb{L}^{5m+5}-9\mathbb{L}^{5m+4}-
-9\mathbb{L}^{5m+3}-6\mathbb{L}^{5m+2}
-3\mathbb{L}^{5m+1}\\
\phantom{[N_4]=}{}+\mathbb{L}^{4m+6}+\mathbb{L}^{4m+5}+2\mathbb{L}^{4m+4}+\mathbb{L}^{4m+3}+\mathbb{L}^{4m+2}\\
\phantom{[N_4]=}{}+2\mathbb{L}^{3m+6}+2\mathbb{L}^{3m+5}+2\mathbb{L}^{3m+4}+2\mathbb{L}^{3m+3}-\mathbb{L}^6.
\end{gather*}
It is left as an exercise to the reader that, for general $d$, a monomial $\mathbb{L}^a$ appears in $[N_d]$ with non-zero coefficient if and only if $a=km+l$ for $d-1\leq k\leq{d\choose 2}$ and ${d\choose 2}-k\leq l\leq{d\choose 2}$.
\item[(5)] The proof of \cite[Lemma~2.6]{BSV} can be read as a motivic identity, resulting in the following formula of a rather different nature:
\[[N_{\bf d}]/[G_{\bf d}]=\sum_{({\bf d}^k)_k}\mathbb{L}^{-\sum_{k<l}\langle{\bf d}^k,{\bf d}^l\rangle-\sum_k\sum_{i\in Q_0}(d^k_i)^2}\prod_k\prod_{i\in Q_0}\frac{\bigl(\mathbb{L}^{-1}\bigr)_{\sum_{\alpha\colon j\rightarrow i}d^k_j}}{\bigl(\mathbb{L}^{-1}\bigr)_{d^k_i}\bigl(\mathbb{L}^{-1}\bigr)_{\sum_{\alpha\colon j\rightarrow i}d^k_j-d^{k+1}_i}},\]
where the sum ranges over all ordered decompositions ${\bf d}=\sum_k{\bf d}^k$.
\end{enumerate}

\section{Hesselink stratification and rational/integral level quiver}\label{s4}

\subsection{Introduction}

We want to calculate the motive of $N_{\mathbf{d}}$ as the corresponding element in the Grothendieck ring of varieties using the Hesselink stratification of the nullcone. To describe the Hesselink stratification of $N_{\mathbf{d}}$, we interpret in detail the results of \cite{LB}, thereby setting up notation. The first new results will only be derived after Proposition \ref{prop LB}.

In fact, we will never use (and consequently, never state) the actual definition of these strata themselves, which is rather indirect (very briefly, one considers the $G_{\bf d}$-saturations of equivalence classes of points in $N_{\bf d}$ under the equivalence relation that their sets of so-called optimal one-parameter subgroups coincide). What is only needed for describing the motive of $N_{\bf d}$ is the isomorphism type of the Hesselink strata, which is described in Theorem~\ref{theo LB}, and an explicit parametrization of the strata.

So we will mainly be concerned with the index set of Hesselink strata $L_{\bf d}$, which we will adapt several times. We will provide an overview over the next steps for more convenience and lucidity. Originally, $L_{\bf d}$ is defined in terms of so-called saturated sets of weights in $R_{\bf d}(Q)$, or, equivalently, in terms of certain so-called dominant minimal coweights $s$. Our first goal is to extend the index set $L_{\mathbf{d}}$ to the set of {\it all} dominant and minimal coweights $\widetilde{L}_{\mathbf{d}}$, see Remark \ref{rem saturated subsets} and the subsequent paragraphs. In a next step we use Theorem~\ref{theo balanced semi-stable empty} to extend $\widetilde{L}_{\mathbf{d}}$ to the set of dominant and balanced coweights $\overline{L}_{\mathbf{d}}$. Afterwards we show that $\overline{L}_{\mathbf{d}}$ is in bijection to $L'_{\mathbf{d}}$, the set of balanced dimension vectors for the so-called rational level quiver $Q_{\mathbb{Q}}$ (see Proposition~\ref{pro par1}) and ultimately we prove that $L'_{\mathbf{d}}$ is in bijection to $L(\mathbf{d})$, a set of tuples of dimension vectors for the so-called integral level quiver, satisfying certain additional requirements (see Proposition~\ref{prop par2}). We summarize and finish this overview with a diagram of index sets
\[L_{\bf d}\subset\widetilde{L}_{\bf d}\subset\overline{L}_{\bf d}\leftrightarrow L_{\bf d}'\leftrightarrow L({\bf d}).\]
Since the Hesselink stratification decomposes the nullcone into finitely many locally closed subsets $S_s$, we obtain the following identity of motives in the Grothendieck ring of varieties:
\begin{equation}\label{iom}
[N_{\mathbf{d}}] = \sum_{s \in L_{\mathbf{d}}} [S_s] .
\end{equation}
Our first goal is to write the set of strata $L_{\bf d}$ in terms of so-called dominant and balanced strings. This description covers all strata, albeit including empty ones. Depending on $s$ we construct a tuple $(Q_s, \theta_s, {\bf d}_s)$ consisting of a~quiver, a stability and a dimension vector, such that the stratum~$S_s$ is non-empty if and only if there exists a $\theta_s$-semistable representations in $R_{{\bf d}_s}(Q_s)$.

We will use and improve the description of the set of strata by Le Bruyn and his construction of $(Q_s, \theta_s, {\bf d}_s)$ \cite{LB}, beginning with the introduction of essential notation:

Starting with a quiver $Q$ with $|Q_0| = n$ vertices, set of arrows $Q_1$, and a dimension vector $\mathbf{d}=(d_i)_i$, we denote $d = \dim{\bf d}$. We decompose the interval
\begin{align*}
 [1 \dots d] = \bigcup_{v=1}^n I_v
\end{align*}
into vertex intervals
\begin{align*}
I_v = \left[ \sum_{i=1}^{v-1} d_i +1 \dots \sum_{i=1}^v d_i\right] .
\end{align*}
Note that we need to fix an enumeration of the vertices of $Q$, but the following constructions will be essentially independent of this choice.

We fix a maximal torus $T_{\mathbf{d}}$ of $G_{\mathbf{d}}$ and denote by $\Pi$ the set of weights of $T_{\mathbf{d}}$ having a nontrivial weight space in $R_{\mathbf{d}}(Q)$. Explicitly, the group of weights of $T_{\mathbf{d}}$ is isomorphic to $\mathbb{Z}^d$ via the canonical generators $\pi_i$ for $1 \leq i \leq d$. Denoting $\pi_{ij} = \pi_j - \pi_i$, we obtain
 \begin{align*}
 \Pi = \{\pi_{ij} = \pi_j - \pi_i \mid i \in I_v,\, j \in I_{v'} \text{ for some arrow } \alpha\colon v \rightarrow v' \text{ in } Q\}
 \end{align*}
 and we denote by $R_{\mathbf{d}}(Q)_{\pi_{ij}}$ the corresponding weight space.
\begin{Definition}\label{subsetR}\quad
\begin{itemize}\itemsep=0pt
\item[(a)] A subset $\R \subset \Pi$ is said to be unstable if there exists a coweight $s=(s_1, \dots, s_d) \in \mathbb{Q}^d$ such that $s_j -s_i \geq 1$ for all $\pi_{ij} \in \R$.
If $\R$ is unstable, there is a unique coweight $s(\R) \in \mathbb{Q}^d$ (see \cite[Lemma 6.2]{LBB}) with this property and such that the norm(-square) $|s(\R)| = s_1^2 + \dots + s_d^2$ is minimal. We will call such a coweight minimal.

\item[(b)] For a coweight $s=(s_1, \dots, s_d) \in \mathbb{Q}^d$ we define the corresponding subset $\R(s) \subset \Pi$ as
\begin{align*}
 \R(s) = \{\pi_{ij} \in \Pi \mid s_j - s_i \geq 1\}.
\end{align*}

\item[(c)] We call $\R$ a saturated subset of $\Pi$ whenever $\R = \R(s(\R))$.

\item[(d)] We denote by $S_{\mathbf{d}} = S_{d_1} \times \dots \times S_{d_n}$ the Weyl group of $G_{\mathbf{d}}$.
\end{itemize}
\end{Definition}

Using this terminology, we can now describe the indexing set for the Hesselink strata; the formal definition will follow in Definition \ref{baldom}.

\begin{Remark}[{see \cite[Section 2.2]{LB}}] \label{rem saturated subsets}
If $\R$ is a saturated subset of $\Pi$, we have a corresponding saturated subspace $X_\R$ of $N_{\mathbf{d}}$ by taking
\begin{align*}
X_\R = \bigoplus_{\pi_{ij} \in \R} R_{\mathbf{d}}(Q)_{\pi_{ij}} .
\end{align*}
By \cite[Section 2.2]{LB}, every Hesselink stratum is a Zariski-open subset of some $X_\R$. We can thus label the Hesselink strata by certain saturated sets of weights $\R$, determined up to conjugacy. Namely, by \cite{Hesse} there is a bijection between the $G_{\mathbf{d}}$-conjugacy classes of saturated subspaces of~$N_{\mathbf{d}}$ and the conjugacy classes of saturated subsets $\R \subset \Pi$ under the action of the Weyl group~$S_{\mathbf{d}}$. 
This correspondence assigns to $\R$ the subspace $X_\R$, and to a saturated subspace the set of nonzero weights of its elements. In particular, the number of Hesselink strata of $N_{\mathbf{d}}$ is at most the number of conjugacy classes of saturated subsets $\R \subset \Pi$ under the Weyl group.

 \end{Remark}
We construct a bijection between the conjugacy classes of saturated subsets $\R \subset \Pi$ under the Weyl group and coweights satisfying certain additional properties.
Consider a coweight $s=(s_1, \dots, s_d) \in \mathbb{Q}^d$ for $G_{\mathbf{d}}$. Then there is a unique partition of $[1 \dots d]$ into a minimal number of pairwise disjoint subsets $\tilde{J}$, such that the $s_j$ for $j\in\tilde{J}$ form an ordered sequence increasing by $0$ or $1$ in each step after sorting. We call such a sequence a segment and denote it by $J$.

A clarifying example: Consider the coweight
\begin{align*}
\left(\frac{3}{2}, 0, -\frac{1}{2}, \frac{4}{3}, 1, - \frac{1}{2}, \frac{1}{3}, \frac{1}{2}, -1,\frac{3}{2}\right) \in \mathbb{Q}^{10} .
\end{align*}
The corresponding segments are
\begin{align*}
J_1 &= \left(-\frac{1}{2}, -\frac{1}{2}, \frac{1}{2}, \frac{3}{2}, \frac{3}{2}\right) , \qquad
J_2 = (-1, 0, 1) ,\qquad
J_3 = \left(\frac{1}{3}, \frac{4}{3}\right)
\end{align*}
and the corresponding partition of $[1 \dots 10]$ is $\tilde{J}_1 = [3, 6, 8, 1, 10]$, $\tilde{J}_2 = [9, 2, 5]$ and $\tilde{J}_3 = [7, 4]$.

After introducing this notation, we can now state definitions for the required conditions on coweights.
\begin{Definition} \label{baldom}\quad
\begin{itemize}\itemsep=0pt
\item[(a)] A coweight is called balanced if and only if for every segment $J$ of $s$ we have
\smash{$
\sum_{i \in \tilde{J}} s_i = 0 $}.
\item[(b)] We call a coweight $s \in \mathbb{Q}^d$ dominant if and only if for every vertex $v$ we have, for all $i,j \in I_v$,
$
s_i \leq s_j $ whenever $ i\leq j $.
\item[(c)] We denote by $L_{\bf d}$ the set of dominant minimal coweights for which $X_{\R(s)}$ is the closure of a Hesselink stratum.

\item[(d)] We denote by $\widetilde{L}_{\mathbf{d}}$ the set of all dominant and minimal coweights.
\end{itemize}
\end{Definition}

Since every saturated subset $\R$ is determined by its unique minimal coweight $s(\R)$, and every minimal coweight belongs to exactly one saturated subset, we have a bijection between the conjugacy classes of saturated subsets $\R \subset \Pi$ under the Weyl group and minimal dominant coweights.

Since the property of being minimal can be hard to verify, our goal is to extend the index set $\widetilde{L}_{\mathbf{d}}$ to $\overline{L}_{\mathbf{d}}$ in such a way, that the new summands in the equation \eqref{iom} are equal to zero, and such that it is easy to verify if a string belongs to the new index set. The property of being balanced will play an essential role in extending the index set of the Hesselink strata. The following lemma portrays the first step \cite{LB}, see also \cite{Hesse}.

\begin{Lemma}
A minimal coweight $\mu = (\mu_1, \dots, \mu_d) \in \mathbb{Q}^d$ is balanced.
\end{Lemma}

\begin{proof}
Suppose $\mu = (\mu_1, \dots, \mu_d)$ is not balanced. Then there exists a segment $J$ such that $\sum_{i \in \tilde{J}} \mu_i \neq 0$. By definition of the segments, for every $\mu_j \in J$ and $\mu_k \notin J$ we have $|\mu_j - \mu_k| < 1$ or $|\mu_j - \mu_k| > 1$. Consequently, we can choose an $ \epsilon >0$ such that these inequalities also hold for~${\mu_j + \epsilon}$ and $\mu_j - \epsilon$. By shifting all entries of $J$ by $\epsilon$ or $-\epsilon$, the distance between two entries in the segment does not change. In this way, we construct another string $\mu'$ such that $\R(\mu) = \R(\mu')$. First assuming $\sum_{i \in \tilde{J}} \mu_i = \alpha > 0$ and we define $\mu'$ by
\begin{gather*}
 \mu_i' = \begin{cases}
 \mu_i & \text{if } i \notin \tilde{J} ,\\
 \mu_i - \epsilon & \text{if } i \in \tilde{J} .
 \end{cases}
\end{gather*}

We denote the cardinality of $\tilde{J}$ by $n$ and obtain
\begin{align*}
|\mu'|- |\mu|&= \sum_{i=1}^d (\mu'_i)^2 - \sum_{i=1}^d (\mu_i)^2 = \sum_{\tilde{J}} (\mu_i - \epsilon)^2 - \sum_{\tilde{J}} (\mu_i)^2 \\
&= \sum_{\tilde{J}} -2\mu_i\epsilon +\epsilon^2 = -2\alpha\epsilon + n \epsilon^2 .
\end{align*}
By choosing $\epsilon < \frac{2\alpha}{n}$, we obtain $|\mu'| < |\mu|$. Accordingly, $\mu$ is not minimal. The case $\sum_{i \in \tilde{J}} \mu_i < 0$ follows analogously with
\begin{gather*}
 \mu_i' = \begin{cases}
 \mu_i & \text{if } i \notin \tilde{J}, \\
 \mu_i + \epsilon & \text{if } i \in \tilde{J} .
 \end{cases}\tag*{\qed}
\end{gather*}\renewcommand{\qed}{}
\end{proof}

We can thus extend the index set of the Hesselink strata using balanced and dominant coweights (this will be done before Lemma \ref{lem balanc stab =0}). The problem of detecting which strata are non-empty remains. To this effect, for every coweight $s$, we define a tuple $(Q_s,{\bf d}_s,\theta_s)$, consisting of a quiver, a dimension vector and a stability, as follows.

\begin{Definition}
For a balanced and dominant coweight $s$ with segments $J_1, \dots, J_u$, we denote by $a_{ik}$ the multiplicity of $p_i +k$ in $J_i$, assuming $J_i$ starts with entry $p_i$. We decompose
\begin{align*}
a_{i,k} = \sum_{v=1}^n a_{i,k}^v ,
\end{align*}
where $a_{i,k}^v$ is the number of entries $a \in I_v$ such that $s_a = p_i +k$.

For each segment $J_i$, where $1 \leq i \leq u$, we consider the level quiver $Q_i$ on $n \times (k_{i} +1)$ vertices~${
\{ (v,j) \mid v \in Q_0 \text{ and } 0\leq j \leq k_i \} }$,
where $p_i +k_i$ is the last entry in $J_i$. In $Q_i$ the number of arrows from $(v,k)$ to $(v',k+1)$ equals the number of arrows from $v$ to $v'$ in $Q$. The dimension vector ${\bf d}_i$ of $Q_i$ is given by $(a_{i,k}^v)_{v,k}$.
\end{Definition}

\begin{Definition}\label{levelquiver}
The quiver $Q_s$ is defined as the disjoint union of the level quivers $Q_i$ associated~to the different segments $J_i$ of $s$ for $1 \leq i \leq u$. The dimension vector ${\bf d}_{s}$ for $Q_{s}$ is the vector~obtained from the dimension vectors ${\bf d}_i$ of $Q_i$.

A stability on $Q_{s}$ is defined by \smash{$\theta_{s} = \bigl(m_{v,j}^{(i)}\bigr)_{i,v,j}$}, where \smash{$m_{v,j}^{(i)} = - z (p_i+j)$} and $z$ is the least common multiple of all denominators of the rational numbers $p_i$ with $1 \leq i \leq u$. In fact, the factor $z$ is to ensure compatibility with \cite{LB}.
\end{Definition}

We recall that a representation $V$ of a quiver $Q$ of dimension vector ${\bf d}$ is called $\theta$-semistable if~${\theta({\bf d})=0}$ and $\theta({\bf e})\leq 0$ whenever ${\bf e}$ is the dimension vector of a subrepresentation of $V$ (note that this definition uses the opposite sign convention to \cite{LB}). We then denote by $R_{\bf d}^{\theta\text{\rm -sst}}(Q)=R_{\bf d}^{\mathrm{sst}}(Q)$ the $\theta$-semistable locus in $R_{\bf d}(Q)$. More generally, using the slope function $\sigma=\theta/{\dim}$ (and not requiring $\theta({\bf d})=0$) we call the representation $V$ $\sigma$-semistable if $\sigma({\bf e})\leq\sigma({\bf d})$ for all dimension vectors ${\bf e}$ of non-zero subrepresentations of $V$.

The following theorem of \cite{LB} summarizes the main properties relating the representation spaces of level quivers and the Hesselink stratification.

\begin{Theorem} \label{theo LB}
The elements $s$ of $\widetilde{L}_{\mathbf{d}}$ give rise to Hesselink strata $S_s$, where $S_s$ is non-empty if and only if the set of $\theta_{s}$-semistable representations in $R_{{\bf d}_s}(Q_s)$ is non-empty.
The Hesselink strata are isomorphic to associated fibre bundles
\smash{$
S_s \cong G_{\mathbf{d}} \times^{P_s} V_s^{\mathrm{sst}}$}.
Here $P_s$ is a certain parabolic subgroup of $G_{\mathbf{d}}$ and $V_s^{\mathrm{sst}}$ is given as follows:

There exists a vector bundle
$
p_s\colon V_s \rightarrow R_{{\bf d}_s} (Q_s)$,
such that
$
V_s^{\mathrm{sst}} = p_s^{-1} \bigl(R_{{\bf d}_s}^{\theta_s\text{\rm -sst}}(Q_s)\bigr)$.
\end{Theorem}

\begin{Proposition} \label{prop LB}
$P_s$, $V_s$ and $R_{{\bf d}_s}(Q_s)$ can be described directly in terms of $s \in \widetilde{L}_{\mathbf{d}}$, namely:
Consider $|Q_1|$-tuples of $d \times d$ matrices. Then $V_s$ is given as the affine space of such tuples, where the matrix corresponding to an arrow $\alpha\colon v\rightarrow v'$ has arbitrary entries at a position $(j,i)$ if~${i\in I_v}$, $j\in I_{v'}$ and $s_j - s_i \geq 1$, and zero entries everywhere else $($compare with Definition {\rm\ref{subsetR})}. The group $P_s$ is given by a single invertible $d \times d$ matrix with zero entries everywhere except possibly at a position $(j,i)$ if $i,j\in I_v$ for some vertex $v\in Q_0$ and $s_j - s_i \geq 0$. For the affine subspace $B_s\subset V_s$ obtained by additionally requiring zero entries at all positions $(j,i)$ except the ones where $s_j - s_i = 1$, we have $B_s \cong R_{{\bf d}_s}(Q_s)$.
\end{Proposition}

Now we can again consider the identity of motives in the localized base ring $\mathbb{M}$ from the beginning of this section
\begin{align*}
[N_{\mathbf{d}}] = \sum_{s \in L_{\mathbf{d}}} [S_s] .
\end{align*}
We divide by the motive of $G_ {\mathbf{d}}$ and use the isomorphism $S_s \cong G_{\mathbf{d}} \times^{P_s} V_s^{\mathrm{sst}}$ to obtain the following:
\begin{align*}
\frac{[ N_{\mathbf{d}}]}{[ G_{\mathbf{d}}]}
&= \sum_{s \in L_{\mathbf{d}}} \frac{[S_s]}{[G_{\mathbf{d}}]}
= \sum_{s \in \widetilde{L}_{\mathbf{d}}} \frac{\big[G_{\mathbf{d}} \times^{P_s} V_s^{\mathrm{sst}}\big]}{[G_{\mathbf{d}}]} .
\end{align*}

Since the group $P_s$ is special, we know that
$
\big[G_{\mathbf{d}} \times^{P_s} V_s^{\mathrm{sst}}\big] \cdot [P_s] = [G_{\mathbf{d}}] \cdot \big[V_s^{\mathrm{sst}}\big]
$
and find
\begin{align*}
\frac{[ N_{\mathbf{d}}]}{[ G_{\mathbf{d}}]} = \sum_{s \in \widetilde{L}_{\mathbf{d}}} \frac{\big[V_s^{\mathrm{sst}}\big]}{[P_s]} .
\end{align*}

Recall from Theorem~\ref{theo LB} that there is a vector bundle
$
p_s \colon V_s \rightarrow R_{\mathbf{d}_s} (Q_s)
$
 of rank $\dim V_s -\dim R_{\mathbf{d}_s}(Q_s)$. Since $R^{\mathrm{sst}}_{\mathbf{d}_s}(Q_s)$ is Zariski-open in $R_{\mathbf{d}_s} (Q_s)$ we also obtain a vector bundle
$
 p'_s\colon V_s^{\mathrm{sst}} \allowbreak\rightarrow R^{\mathrm{sst}}_{\mathbf{d}_s}(Q_s)
$
 of the same rank as long as $ R^{\mathrm{sst}}_{\mathbf{d}_s}(Q_s)$ is not empty. Moreover, the reductive part of~$P_s$ is isomorphic to $G_{\mathbf{d}_s}$, yielding
\begin{align*}
\frac{[ N_{\mathbf{d}}]}{[ G_{\mathbf{d}}]} =\sum_{s \in \widetilde{L}_{\mathbf{d}}} \frac{\big[R^{\mathrm{sst}}_{\mathbf{d}_s}(Q_s)\big] \cdot \mathbb{L}^{\dim V_s -\dim R_{\mathbf{d}_s}(Q_s)}}{[G_{\mathbf{d}_s}] \cdot \mathbb{L}^{\dim P_s - \dim G_{\mathbf{d}_s}}}.
\end{align*}
We denote $\phi(s) = \dim V_s -\dim R_{\mathbf{d}_s}(Q_s) -( \dim P_s - \dim G_{\mathbf{d}_s})$ and simplify the equation to
\begin{align} \label{formula motive}
\frac{[ N_{\mathbf{d}}]}{[ G_{\mathbf{d}}]} = \sum_{s \in {\widetilde{L}_{\mathbf{d}}}} \mathbb{L}^{\phi(s)} \frac{\big[R^{\mathrm{sst}}_{\mathbf{d}_s}(Q_s)\big]}{[G_{\mathbf{d}_s}]} .
\end{align}

Now we extend the index set $\widetilde{L}_{\mathbf{d}}$ of dominant and minimal strings to the index set $\overline{L}_{\mathbf{d}}$ of dominant and balanced strings. Theorem~\ref{theo balanced semi-stable empty} below shows that this change of index sets does not affect the formula for the motive of the nullcone. We need some preparation in advance.

\begin{Lemma} \label{lem balanc stab =0}
Let $\mu$ be a balanced string and $J_i$ a segment of $\mu$ starting with $p_i$. We consider the stability $\tilde{\theta}$ defined by $\tilde{\theta}(v,j) = -(p_i+j)$. Then for the dimension vector $\mathbf{d}_{J_i}$ we have $\tilde{\theta} (\mathbf{d}_{J_i})=0$.
\end{Lemma}

\begin{proof}We have
\begin{align*}
\tilde{\theta} (\mathbf{d}_{J_i})= \sum_{(v,j) \in Q_{J_i}} - a_{i,k}^v (p_j+j) = \sum_{j=1}^{k_i} a_{i,k} (p_j+j) = \sum_{k \in \tilde{J}_i} \mu_k =0 ,
\end{align*}
because $\mu$ is balanced.
\end{proof}

\begin{Definition}For a segment $K$, we define $Q_K=Q_i$ and ${\bf d}_K={\bf d}_i$ if $K=J_i$.

Let $\mu$ be a string and $J$ a segment of $\mu$ with corresponding index set $I$. For a subset $\hat{I}$ of~$I$, we can construct a subquiver $\hat{Q}_J$ of $Q_J$ from the substring $\hat{J}$ of $J$ determined by $\hat{I}$ in the same way the level quiver is constructed from a balanced string; we can also construct a dimension vector $\hat{{\bf d}}_J$ for this subquiver which fulfills \smash{$\bigl(\hat{{\bf d}}_J\bigr)_i \leq ({\bf d}_J)_i$}.

We call a subset $\hat{I}$ of $I$ left-complete if, for every $k \in \hat{I}$, all $l$ with $\mu_k -\mu_l =1$ and $\pi_{lk} \in \Pi$ also belong to $\hat{I}$. Note that this implies the existence of an arrow from the vertex corresponding to $\mu_l$ to the vertex corresponding to $\mu_k$ in $Q_J$.
\end{Definition}

\begin{Lemma} \label{lem inequ for left-complete}
If the semistable locus \smash{$R_{d_J}^{\tilde{\theta}\text{\rm -sst}} (Q_J)$} is non-empty, then for all left-complete subsets~$\hat{I}$ of the index set~$I$ corresponding to the segment~$J$ the inequality
$\sum_{\hat{I}} \mu_i \leq 0$
holds.
\end{Lemma}
\begin{proof}Assume, to the contrary, that \smash{$\sum_{\hat{I}} \mu_i > 0$}. Let $V$ be an arbitrary representation of dimension vector ${\bf d}_J$ of $Q_J$. Since $\hat{I}$ is left-complete,we have \smash{$\bigl({\bf d}_J -\hat{\bf d}_J\bigr)_k=0$} for all vertices~$v_k$ in $\hat{Q}_J$, which are not a sink in $\hat{Q}_J$, and \smash{$\bigl({\bf d}_J -\hat{\bf d}_J\bigr)_k=({\bf d}_J)_k$} for all $v_k$ which are not in $\hat{Q}_J$. Therefore, $V$ admits a subrepresentation $V'$ of dimension vector ${\bf d}_J -\hat{\bf d}_J$. By Lemma \ref{lem balanc stab =0}, we have $\tilde{\theta} (\mathbf{d}_J) = 0$ since $\mu$ is balanced, showing that \smash{$\tilde{\theta} \bigl(\mathbf{d}_J - \hat{\bf d}_J\bigr) >0$}. This proves the claim.
\end{proof}

\begin{Theorem} \label{theo balanced semi-stable empty}
For a balanced string $s$ which is not minimal, the semistable locus \smash{$R_{{\bf d}_s}^{\theta_s\text{\rm -sst}}(Q_s)$} is empty.
\end{Theorem}

\begin{proof}
First note that we can switch from the stability $\theta_s$ to the stability $\tilde{\theta}$ given by $\tilde{\theta}(v,j) = -(p_i +j)$ on the level quiver $Q_K$ for a segment $K$ starting with $p_i$, since the open subset of semistable representations does not change when passing from $\theta$ to $x \cdot \theta$ for a positive $x$.

Moreover, we note that for a fixed string $s = (s_1, \dots, s_n) \in \mathbb{Q}^n$, the minimum of $\big\{|s-q| := \sum_{i=1}^n (s_i - q)^2 \big\}$ for $q \in \mathbb{Q}$ is given by $|s-\bar{q}|$ such that $\sum_{i=1}^n (s_i-\bar{q}) =0$. If, additionally, $\sum_{i=1}^n s_i < \sum_{i=1}^n (s_i-q) <0$, then $|s| \geq |s-q|$, since we can consider $|s-q|$ (for a fixed $s$) as a~quadratic function in $q$. Analogously, if $\sum_{i=1}^n s_i > \sum_{i=1}^n (s_i-q) >0$, then $|s| \geq |s-q|$.

Let $\mu = (\mu_1, \dots, \mu_d)$ be a balanced, but not minimal string. By definition, there exists a~(balanced) $\mu'$ with $|\mu'| < |\mu|$ and $\R(\mu') \supset \R(\mu)$. Therefore, we find a segment $J= (\mu_i)_{i \in I}$ in $\mu$ such that the corresponding substring $J' = (\mu')_{i \in I}$ of $\mu'$ satisfies $|J'|= \sum_I \mu_i'^2 < \sum_I \mu_i^2=|J|$.
We can decompose the indexing set $I$ in three disjoint sets given by
\begin{align*}
I^n_1 = \{ i \in I \text{ with } \mu_i > \mu_i' \} ,\qquad
I^p_1 = \{ i \in I \text{ with } \mu_i <\mu_i' \}, \qquad
I_0 = \{ i \in I \text{ with } \mu_i = \mu_i' \} .
\end{align*}
Since $|J'| < |J|$, either $\sum_{I^n_1} \mu_i'^2 < \sum_{I^n_1} \mu_i^2$ or $\sum_{I^p_1} \mu_i'^2 < \sum_{I^p_1} \mu_i^2$. Without loss of generality, we assume the inequality holds for $I^n_1$. The other case follows analogously.

Consider $\mu_i'$ with $i \in I^n_1$. Then there exists an $\epsilon_i >0$ with $\mu_i' = \mu_i - \epsilon_i$. In the level quiver $Q_J$, there exists a vertex $v_i$ corresponding to $\mu_i$. Every arrow $\alpha$ in $Q_J$ with $t(\alpha)= v_i$ corresponds to (possibly several) elements $\mu_l \in J$ such that $\mu_i - \mu_l = 1$ and $\pi_{li} \in \Pi$. Since $\R(\mu') \supset \R(\mu)$, every $\mu_l$ obtained in the above way fulfills $\mu'_l = \mu_l - \epsilon_l$ with $\epsilon_l > \epsilon_i$. This shows that $I^n_1$ is left-complete.

We now assume that the semistable locus in question is not empty.
By Lemma \ref{lem inequ for left-complete}, we know that for all left-complete subset $\hat{I}$ the inequality \smash{$\sum_{\hat{I}} \mu_i \leq 0$} hold.
We know that $\smash{\sum_{I^n_1} \mu_i'^2} = \sum_{I^n_1} (\mu_i - \epsilon_i)^2 < \sum_{I^n_1} \mu_i^2$ and $\sum_{I^n_1} \mu_i \leq 0$, because $I^n_1$ is left-complete. We choose $\tilde{\epsilon}_1$ as the minimum of all $\epsilon_i$ for $i \in I^n_1$. We get
\begin{align*}
\sum_{I^n_1} \mu_i -\epsilon_i < \sum_{I^n_1} \mu_i - \epsilon_i + \tilde{\epsilon}_1 \leq \sum_{I^n_1} \mu_i \leq 0 ,
\end{align*}
leading to \smash{$\sum_{I^n_1} (\mu_i - \epsilon_i)^2 > \sum_{I^n_1} (\mu_i - \epsilon_i + \tilde{\epsilon}_1)^2$}.
Now we define
$
I^n_2 = \{i \in I^n_1 \mid \mu_i -\epsilon_i + \tilde{\epsilon}_1 < \mu_i\} $.
In other words, we exclude all $\mu_i'$ with $\mu_i' + \tilde{\epsilon}_1 = \mu_i$.
Now $I^n_2$ is left-complete by the same argument as for $I^n_1$. We choose $\tilde{\epsilon}_2$ as the minimum of all $\epsilon_i - \tilde{\epsilon}_1$ with $i \in I^n_2$. We get
\begin{align*}
\sum_{I^n_2} \mu_i - \epsilon_i +\tilde{\epsilon}_1 < \sum_{I^n_2} \mu_i - \epsilon_i+ \tilde{\epsilon}_1 + \tilde{\epsilon}_2 \leq \sum_{I^n_2} \mu_i \leq 0
\end{align*}
concluding
\begin{align*}
\sum_{I^n_2} (\mu_i - \epsilon_i + \tilde{\epsilon}_1)^2 > \sum_{I^n_2} (\mu_i - \epsilon_i + \tilde{\epsilon}_1 + \tilde{\epsilon}_2)^2.
\end{align*}
We get a descending chain of subsets
$
I_1^n \supsetneq I_2^n \supsetneq \cdots$,
which has to terminate after finitely many steps. This yields the following estimate:
\begin{align*}
\sum_{I^n_1} (\mu_i -\epsilon_i)^2 > \sum_{I^n_1} (\mu_i - \epsilon_i + \tilde{\epsilon}_1)^2 > \sum_{I^n_1 \setminus I^n_2} \mu_i^2 + \sum_{I^n_2} (\mu_i - \epsilon_i + \tilde{\epsilon}_1 + \tilde{\epsilon}_2)^2 > \dots > \sum_{I^n_1} \mu_i^2 .
\end{align*}
Since we have started with $\sum_{I^n_1} (\mu_i -\epsilon_i)^2 < \sum_{I^n_1} \mu_i^2$, we have thus constructed a contradiction, proving that a balanced, but not minimal string, has an empty semistable locus.
\end{proof}

\begin{Corollary}
 For the nullcone $N_{\bf d}$, we have the following equality of motives:
 \begin{align*}
 \frac{[ N_{\mathbf{d}}]}{[ G_{\mathbf{d}}]} = \sum_{s \in {\overline{L}_{\mathbf{d}}}} \mathbb{L}^{\phi(s)} \frac{\big[R^{\mathrm{sst}}_{\mathbf{d}_s}(Q_s)\big]}{[G_{\mathbf{d}_s}]} .
 \end{align*}
\end{Corollary}

\begin{Example}
We illustrate these constructions and the calculation of the motive of $N_{d}$ for the $m$-loop quiver for $d=2$ and $d=3$.
For $d=2$, the list of strings consists of
$\bigl(-\frac{1}{2}, \frac{1}{2}\bigr)$ and~$ (0,0)$. The corresponding quivers $Q_s$ are
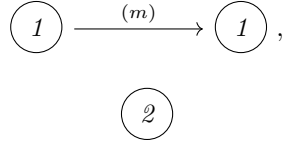
\[\begin{tikzcd}
	\circlearound{1} && \circlearound{1}
	\arrow[from=1-1, to=1-3, "(m)"],
\end{tikzcd}\]
\[\begin{tikzcd}
	\circlearound{2}
\end{tikzcd}\]
with stability $\theta_s = (1,-1)$ and $\theta_s = (0)$, respectively.

We calculate the individual summands in formula \eqref{formula motive} starting with the string $s=\bigl(-\frac{1}{2}, \frac{1}{2}\bigr)$ and thus ${\bf d}_s = (1,1)$:
the space $V_s$ consists of $m$-tuples of matrices with one arbitrary entry, consequently $\dim V_s = m$. Moreover, the dimension of the representation space $R_{{\bf d}_s} (Q_s)$ equals the number of arrows. The group $P_s$ consists of an invertible $2\times 2$ matrix with arbitrary entries at $3$ positions and $\dim G_{d_s} = \dim (G_1 \times G_1) = 2$.
We thus find
\begin{align*}
\phi\left(\left(-\frac{1}{2}, \frac{1}{2}\right)\right) = m - m -3+2 = -1 ,
\end{align*}
and similarly
$
\phi ((0,0)) = 0-0-4+4 = 0 $.
For the first quiver, every representation is semistable except the one consisting only of trivial maps. For the second quiver, there is only one representation, which is in addition semistable. Therefore, the corresponding motives of semistable loci are $\mathbb{L}^m -1$ and $1$, respectively. Consequently, the motive of $N_2$ calculates to
\begin{align*}
[N_2] &= \left(\mathbb{L}^{-1} \frac{(\mathbb{L}^m -1)}{[G_1]^2}+ \frac{1}{[G_2]} \right) [G_2] = \bigl(\mathbb{L}^{m-1} - \mathbb{L}^{-1}\bigr) \bigl(\mathbb{L}^2 + \mathbb{L}\bigr) +1 = \mathbb{L}^{m+1} + \mathbb{L}^m - \mathbb{L} .
\end{align*}

The next example is the $m$-loop quiver with dimension vector $d =3$. The relevant strings in the classification are
\begin{align*}
(-1,0,1), \left(-\frac{2}{3},\frac{1}{3},\frac{1}{3}\right)\left(-\frac{1}{3},-\frac{1}{3},\frac{2}{3}\right) \left(-\frac{1}{2},0, \frac{1}{2}\right), (0,0,0)
\end{align*}
with corresponding quivers
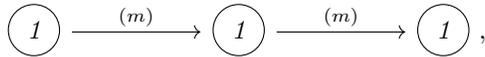
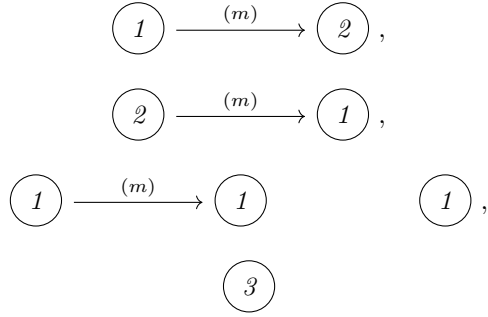
\[\begin{tikzcd}
	\circlearound{1} && \circlearound{1} && \circlearound{1}
	\arrow[from=1-1, to=1-3, "(m)"] \arrow[from=1-3, to=1-5, "(m)"] ,
\end{tikzcd}\]
\[\begin{tikzcd}
	\circlearound{1} && \circlearound{2}
	\arrow[from=1-1, to=1-3, "(m)"] ,
\end{tikzcd}\]
\[\begin{tikzcd}
	\circlearound{2} && \circlearound{1}
	\arrow[from=1-1, to=1-3, "(m)"] ,
\end{tikzcd}\]
\[\begin{tikzcd}
	\circlearound{1} && \circlearound{1} && \circlearound{1}
	\arrow[from=1-1, to=1-3, "(m)"] ,
\end{tikzcd}\]
\[\begin{tikzcd}
	\circlearound{3}
\end{tikzcd}\]
and stability
$(1,0,-1)$, $(2,-1)$, $(1,-2)$, $(1,-1,0)$, $(0) $.
Similarly to the previous example, we find
\begin{gather*}
\phi((-1,0,1)) = 3m -2m -6+3 = m -3 ,\\
\phi\left(\left(-\frac{2}{3}\right), \frac{1}{3},\frac{1}{3}\right) = 2m -2m -7+5 = -2, \\
\phi\left(\left(-\frac{1}{3},-\frac{1}{3},\frac{2}{3},\right)\right) = 2m -2m - 7+ 5 = -2,\\
\phi\left(\left(-\frac{1}{2},0, \frac{1}{2}\right)\right) = m-m -6+3 = -3,\qquad
\phi((0,0,0)) = 0 -0 -9+9 = 0.
\end{gather*}

It remains to describe the motives of the semistable loci. For the first quiver, all representations are semistable unless one tuple of parallel arrows contains only trivial maps. Consequently, the resulting motive is $(\mathbb{L}^m -1)^2$. For the second and third quiver, we can write the $m$-tuple of vectors of length $2$ in a $2 \times m$ matrix, obtaining that the representation is semistable if and only if this matrix is of rank $2$. Consequently, the motive equals $(\mathbb{L}^m -1 )(\mathbb{L}^m - \mathbb{L})$. The motives for the final two quivers are $\mathbb{L}^m -1$ and $1$, respectively, with similar arguments as before.

Putting all this information together, the motive of $N_3$ is given by
\begin{align*}
[N_3]={}& \left( \mathbb{L}^{m-3} \frac{(\mathbb{L}^{m}-1)^2}{[G_1]^3} + 2 \mathbb{L}^{-2} \frac{(\mathbb{L}^m -1)(\mathbb{L}^m-\mathbb{L})}{[G_1][G_2]} + \mathbb{L}^{-3} \frac{\mathbb{L}^m -1}{[G_1]^3} + \frac{1}{[G_3]}\right) [G_3] \\
={}& \bigl(\mathbb{L}^{3m-3} -2 \mathbb{L}^{2m-3}+\mathbb{L}^{m-3}\bigr) \bigl( \mathbb{L}^6 + 2 \mathbb{L}^5 +2 \mathbb{L}^4 + \mathbb{L}^3\bigr) \\
& + \bigl(2\mathbb{L}^{2m-2} - 2 \mathbb{L}^{m-1}-2 \mathbb{L}^{m-2}+2\mathbb{L}^{-1}\bigr) \bigl( \mathbb{L}^4 + \mathbb{L}^3 + \mathbb{L}^2 \bigr) \\
& + \bigl(\mathbb{L}^{m-3} - \mathbb{L}^{-3}\bigr) \bigl(\mathbb{L}^6 + 2 \mathbb{L}^5 +2 \mathbb{L}^4 + \mathbb{L}^3 \bigr) +1\\
={}& \mathbb{L}^{3m+3} + 2\mathbb{L}^{3m+2} + 2\mathbb{L}^{3m+1} + \mathbb{L}^{3m}-2\mathbb{L}^{2m+3} - 4\mathbb{L}^{2m+2} - 4\mathbb{L}^{2m+1} - 2\mathbb{L}^{2m} \\
& + \mathbb{L}^{m+3} + 2\mathbb{L}^{m+2} +2\mathbb{L}^{m+1} + \mathbb{L}^{m} \\
& + 2 \mathbb{L}^{2m+2} +2 \mathbb{L}^{2m+1} +2 \mathbb{L}^{2m} -
2 \mathbb{L}^{m+3} -2 \mathbb{L}^{m+2} -2 \mathbb{L}^{m+1} \\
& - 2 \mathbb{L}^{m+2} -2 \mathbb{L}^{m+1} -2 \mathbb{L}^{m} + 2 \mathbb{L}^{3} + 2 \mathbb{L}^2 +2 \mathbb{L} \\
& + \mathbb{L}^{m+3} + 2 \mathbb{L}^{m+2} + 2\mathbb{L}^{m+1} + \mathbb{L}^{m} - \mathbb{L}^3 - 2 \mathbb{L}^2 - 2 \mathbb{L} - 1 +1 \\
={}& \mathbb{L}^{3m+3} + 2\mathbb{L}^{3m+2} + 2\mathbb{L}^{3m+1} + \mathbb{L}^{3m} - 2\mathbb{L}^{2m+3} - 2\mathbb{L}^{2m+2} -2 \mathbb{L}^{2m+1} + \mathbb{L}^3 ,
\end{align*}
 in accordance with the example in Section \ref{ex part 1}.
\end{Example}

\subsection{Rational/integral level quiver}

We know that the elements $s$ of the index set ${\overline{L}_{\mathbf{d}}}$ give rise to quivers $Q_s$, such that the corresponding Hesselink stratum is non-empty if and only if the semi-stable locus is non-empty. Our next aim is now to define a so-called rational level quiver which combines the level quivers associated to the segments, in such a way that the previous properties are still preserved, but the index set for the Hesselink stratification can be described in terms of specific dimension vectors~$\mathbf{e}$ for this fixed quiver.

\begin{Definition}
For a quiver $Q$, we define the rational level quiver $Q_{\mathbb{Q}}$ with set of vertices~${Q_0 \times \mathbb{Q}}$ and arrows
$
(\alpha, a) \colon (i,a) \rightarrow (j, a+1)
$
for every arrow~${\alpha\colon i \rightarrow j}$ in $Q$ and every~${a \in \mathbb{Q}}$. This quiver obviously admits a decomposition into components $Q_{\mu}$ for $\mu \in [0,1] \cap \mathbb{Q}$, where $Q_{\mu}$ is the full subquiver of all vertices $(i,a)$ such that the fractional part $\{a\}$ of $a$ equals~$\mu$.

We define a stability on $Q_{\mathbb{Q}}$ by $\theta(i,a) = -a$.
We call a dimension vector $\mathbf{e}$ for $Q_\mathbb{Q}$ balanced if $\theta(\mathbf{e}|_{Q_\mu})=0$
for all $\mu\in[0,1]\cap\mathbb{Q}$ (where $\mathbf{e}|_{Q_\mu}$ denotes the restriction of $\mathbf{e}$ to the component~$Q_\mu$). Every finitely supported dimension vector $\mathbf{e}$ for $Q_\mathbb{Q}$ admits a projection $|\mathbf{e}|$ to a dimension vector for $Q$ by $|\mathbf{e}|_i=\sum_{a\in\mathbb{Q}}e_{i,a}.$
\end{Definition}

\begin{Proposition}\label{pro par1}
Using this notation, the index set ${\overline{L}_{\mathbf{d}}}$ for the Hesselink stratification of $N_{\mathbf{d}}$ is in bijection to the set of all balanced dimension vectors $\mathbf{e}$ for $Q_\mathbb{Q}$ such that $|\mathbf{e}|=\mathbf{d}$. Moreover, we have an isomorphism of semistable loci \smash{$R^{\theta_s\text{\rm -sst}}_{\mathbf{d}_s}(Q_s) \cong R^{\theta\text{\rm -sst}}_{\mathbf{e}}(Q_{\mathbb{Q}})$}.
\end{Proposition}

\begin{proof}
We construct a dominant coweight from a balanced dimension vector $\mathbf{e}$ in the following way:

For every entry $e_{i,a}$ in $\mathbf{e}$, we take $e_{i,a}$-times the entry $a$ in the coweight, in weakly ascending order for every vertex of $Q$. That this induces a bijection between balanced dimension vectors and balanced and dominant coweights now follows immediately from the definitions.

Now we have to compare the stabilities given by $\theta_s$ and $\theta$. Relabelling the vertices $(v,j)$ of the quiver $Q_s$ to $(v, p_i+j)$, we see directly that $Q_s$ can be viewed as a subquiver of the rational level quiver. Using the above construction of a balanced dominant coweight from a~balanced dimension vector $\mathbf{e}$, it is then clear that in both cases the multiplicity of the entry~${p_i +j}$ (resp.~$a$) in the string is encoded in the dimension vector. So the only essential difference between the stabilities is the factor $z$ of Definition \ref{levelquiver}. Since the set of semistable representations with respect to a stability $\theta$ does not change when passing from $\theta$ to $x \cdot \theta$ for positive $x$, it follows that, in fact, $R^{\theta_s\text{\rm -sst}}_{\mathbf{d}_s}(Q_s) \cong R^{\theta\text{\rm -sst}}_{\mathbf{e}}(Q_{\mathbb{Q}})$.
\end{proof}

We denote the set of balanced dimension vectors $\mathbf{e}$ for ${Q}_{\mathbb{Q}}$ such that $|{\bf e}| = {\bf d}$ by $L'_\mathbf{d}$.

The quiver $Q_\mathbb{Q}$ is the disjoint union of the $Q_\mu$ for $\mu\in[0,1]\cap\mathbb{Q}$, and every $Q_\mu$ can be identified with $Q_{\mu=0}=Q_\mathbb{Z}$, called the integral level quiver. We can reparametrize the index set $L'_{\bf d}$: we map $\mathbf{e}$ to the tuple $(\mathbf{e}|_{Q_\mu})_\mu$ of dimension vectors, which we view as dimension vectors for $Q_\mathbb{Z}$. 
Now we interpret this reparametrization in terms of our stability conditions:

After identifying a balanced dimension vector $\mathbf{e}$ for $Q_{\mathbb{Q}}$ with a tuple of dimension vectors~$(\mathbf{e}_\mu)_\mu$ for $Q_{\mathbb{Z}}$, we have
$
\theta (\mathbf{e}|_{Q_\mu}) = \theta (\mathbf{e}_\mu) -\mu\cdot\dim \mathbf{e}_\mu $.
This condition can be rephrased using the slope function $\sigma=\theta/{\dim}$ and noting that $\theta (\mathbf{e}|_{\mathbb{Q}_{\mu}}) = 0$ since $\mathbf{e}$ is balanced
\begin{align*}
\sigma ( \mathbf{e}_\mu) = \frac{\theta(\mathbf{e}_{\mu})}{\dim \mathbf{e}_\mu} = \mu
\end{align*}
if ${\bf e}_\mu\not=0$.

We can now state a final reformulation for the index set of Hesselink strata.

\begin{Proposition} \label{prop par2}
The index set $L'_{\bf d}$ for the Hesselink stratification can be viewed as the set of tuples $(\mathbf{e}_\mu)_\mu$ of dimension vectors for $Q_\mathbb{Z}$ such that ${\bf e}_\mu=0$ or $\sigma(\mathbf{e}_\mu)
=\mu$ for all $\mu$, and such that~${\sum_\mu|\mathbf{e}_\mu|=\mathbf{d}}$. We will denote the set of such tuples by $L({\bf d})$ in the following.
\end{Proposition}

\begin{proof}
The index sets $L'_{\bf d}$ and $L({\bf d})$ are obviously in bijection by the previous construction. Similarly to the proof of Proposition~\ref{pro par1}, the semistable locus does not change when passing from $\theta$ to $\theta - y\cdot\dim $ for an arbitrary $y$.
\end{proof}

\section{Wall-crossing type formula}\label{s5}

Here we understand the notion of a wall-crossing formula in a purely formal sense: an identity of generating functions of motives involving an ordered product over slopes. A prototypical example is the identity \cite[Theorem 2.4]{MozIntro} expressing the generating series of motives of moduli stacks of quiver representations as an ordered product over such series for semistable quiver representations of fixed slope.

Our aim in this section is thus to interpret the sum formula \eqref{formula motive} as an ordered product formula for the corresponding generating series.

Describing the strata in terms of balanced dimension vectors $\mathbf{e}$ for the quiver $Q_\mathbb{Q}$, we can rewrite the formula
$
\phi(s) = \dim V_s - \dim R_{\mathbf{d}_s} (Q_s) - (\dim P_s - \dim G_{\mathbf{d}_s})
$
in terms of $\mathbf{e}$.

Using Proposition \ref{prop LB}, we see that $ \dim V_s - \dim R_{\mathbf{d}_s} (Q_s)$ is the dimension of the space of $|Q_1|$-tuples of $d\times d$ matrices with zero entries except possibly at position $(l,k)$ in the matrix corresponding to $\alpha\colon i\rightarrow j$ if $k\in I_i$, $l\in I_{j}$ and $s_l - s_k > 1$. In terms of the rational level quiver, this dimension is encoded in the dimensions $e_{i,a}$ at vertices which are not directly adjacent
\begin{align*}
\dim V_s - \dim R_{\mathbf{d}_s} (Q_s) = \sum_{\alpha\colon i \rightarrow j} \sum_{a+1<b} e_{i,a} e_{j,b} .
\end{align*}
We analyse the summand $\dim P_s - \dim G_{\mathbf{d}_s}$. The group $P_s$ was described in terms of matrix entries $(l,k)$ such that $k,l\in I_i$ for some vertex $i\in Q_0$ and $s_l - s_k \geq 0$, while $G_{\mathbf{d}_s}$ is determined by $s_l - s_k = 0$. So we can construct $P_s$ from $G_{\mathbf{d}_s}$ by allowing arbitrary entries at positions $(l,k)$ where $k,l\in I_i$ for some vertex $i\in Q_0$ and $s_l - s_k > 0$. For symmetry reasons, the number of such positions is exactly
\begin{align*}
\frac {\sum_{i\in Q_0} d_i^2 - \dim G_{\mathbf{d}_s}}{2} .
\end{align*}
\begin{Remark}
Alternatively, we can note that $\dim P_s - \dim G_{\mathbf{d}_s}$ equals the dimension of the unipotent radical corresponding to the parabolic subgroup $P_s$.
\end{Remark}
With \smash{$d_i = \sum_{a \in \mathbb{Q}} e_{i,a}$} and \smash{$ \dim G_{\mathbf{d}_s} = \sum_{ i \in Q_0} \sum_{a \in \mathbb{Q}} e_{i,a}^2$}, we can conclude
\begin{align*}
\dim P_s - \dim G_{\mathbf{d}_s} = \frac{\sum_{ i \in Q_0} \bigl(\sum_{a \in \mathbb{Q}} e_{i,a}\bigr)^2 - \sum_{ i \in Q_0} \sum_{a \in \mathbb{Q}} e_{i,a}^2 } {2} = \sum_{i \in Q_0} \sum_{a < b} e_{i,a} e_{i,b}
\end{align*}
yielding a formula for $\phi$ in terms of the dimension vector $\mathbf{e}$
\begin{align*}
\phi(\mathbf{e}) = \sum_{\alpha\colon i \rightarrow j} \sum_{a+1 <b} e_{i,a} e_{j,b} - \sum_{ i \in Q_0} \sum_{a < b} e_{i,a} e_{i,b} .
\end{align*}
Since we have shifted the dimension vector $\mathbf{e}$ to a tuple $(\mathbf{e}_{\mu})_\mu$ of dimension vectors for $Q_\mathbb{Z}$, we can rewrite this formula in terms of $(\mathbf{e}_{\mu})_{\mu}$. In the first summand, we replace $e_{i,a}$ and $e_{j,b}$ by~\smash{$e_{i,a-\mu}^{\mu}$} and \smash{$e_{j,b-\mu'}^{\mu'}$}, respectively, and similarly in the second term. Note that, consequently, $a$ and $b$ denote integers in all the following summations:
\begin{align*}
\phi((\mathbf{e}_\mu)_\mu) &= \sum_{\alpha\colon i \rightarrow j} \sum_{a+1 <b}\sum_{\mu,\mu'} e^\mu_{i,a-\mu} e^{\mu'}_{j,b-\mu'} - \sum_{ i \in Q_0} \sum_{a < b}\sum_{\mu,\mu'} e^\mu_{i,a-\mu} e^{\mu'}_{i,b-\mu'} \\
&= \sum_{\alpha\colon i \rightarrow j} \sum_{a+1+\mu <b+\mu'}\sum_{\mu,\mu'} e^\mu_{i,a} e^{\mu'}_{j,b} - \sum_{ i \in Q_0} \sum_{a+\mu < b+\mu'}\sum_{\mu,\mu'} e^\mu_{i,a} e^{\mu'}_{i,b}.
\end{align*}

The condition $a+1+\mu < b+\mu' $ is equivalent to ($\mu < \mu'$ and $ a+1=b$) or $a+1<b$ for~${\mu, \mu' \in [0,1]}$ and $a,b \in \mathbb{Z}$. Similarly, $a+\mu<b+\mu'$ is equivalent to ($\mu < \mu'$ and $ a=b$) or~${a<b}$.
We thus obtain
\begin{align*}
\phi((\mathbf{e}_{\mu})_{\mu}) ={}& \sum_{\alpha\colon i \rightarrow j} \biggl(\sum_{a+1 < b} \sum_{\mu, \mu'} e_{i,a}^{\mu} e_{j,b}^{\mu'} + \sum_{a+1 = b} \sum_{\mu < \mu'} e_{i,a}^{\mu} e_{j,b}^{\mu'}\biggr) \\
&- \sum_{i} \biggl(\sum_{a < b} \sum_{\mu, \mu'} e_{i,a}^{\mu} e_{i,b}^{\mu'} + \sum_{a = b} \sum_{\mu < \mu'} e_{i,a}^{\mu} e_{i,b}^{\mu'}\biggr) .
\end{align*}

With our new parametrization of the index set and Propositions \ref{pro par1} and \ref{prop par2}, we can now state the following formula for the motive of $N_{\mathbf{d}}$.

\begin{Theorem} \label{theo mot}
For the nullcone $N_{\mathbf{d}}$ of a quiver $Q$ with dimension vector $\mathbf{d}$, we have the following identities of motives:
\begin{align*}
\begin{split}
\frac{[N_{\mathbf{d}}]}{[G_{\mathbf{d}}]}
&= \sum_{s \in \overline{L}_{\mathbf{d}}} \mathbb{L}^{\phi(s)} \frac{\big[R^{\theta_s\text{\rm -sst}}_{\mathbf{d}_s}(Q_s)\big]}{[G_{\mathbf{d}_s}]} = \sum_{\mathbf{e} \in L'_{\mathbf{d}}} \mathbb{L}^{\phi(\mathbf{e})} \frac{\big[R_{\mathbf{e}}^{\theta\text{\rm -sst}}(Q_{\mathbb{Q}})\big]}{[G_{\mathbf{e}}]} \\
&= \sum_{(\mathbf{e}_\mu)_\mu\in L(\mathbf{d})} \mathbb{L}^{\phi((\mathbf{e}_\mu)_\mu)}\prod_{\mu} \frac{\big[R_{\mathbf{e}_{\mu}}^{\theta\text{\rm -sst}} (Q_{\mathbb{Z}})\big]}{[G_{\mathbf{e}_{\mu}}]} .
\end{split}
\end{align*}
\end{Theorem}

\begin{Definition}
We define the $k$-shift ${\tau}^k \mathbf{e}$ of a dimension vector $\mathbf{e}$ for $Q_\mathbb{Z}$ as the dimension vector
$
\bigl( {\tau}^k \mathbf{e}\bigr)_{i,a} = \mathbf{e}_{i,a+k} $.
\end{Definition}

A brief calculation shows that for the Euler form of $Q_\mathbb{Z}$, we have
\begin{align*}
-\big\langle\mathbf{e}^{\mu}, {\tau}^k \mathbf{e}^{\mu'}\big\rangle = - \sum_i \sum_a e_{i,a}^{\mu} e_{i, a+k}^{\mu'} + \sum_{i \rightarrow j} \sum_a e_{i,a}^{\mu} e_{j,a+1+k}^{\mu'} .
\end{align*}

We now consider $\phi((\mathbf{e}_\mu)_\mu)$, which equals
\begin{gather*}
\sum_{\alpha\colon i \rightarrow j} \biggl(\sum_{a+1 < b} \sum_{\mu, \mu'} e_{i,a}^{\mu} e_{j,b}^{\mu'} + \sum_{a+1 = b} \sum_{\mu < \mu'} e_{i,a}^{\mu} e_{j,b}^{\mu'}\biggr)\\
\qquad-\sum_{i} \biggl(\sum_{a < b} \sum_{\mu, \mu'} e_{i,a}^{\mu} e_{i,b}^{\mu'} + \sum_{a = b} \sum_{\mu < \mu'} e_{i,a}^{\mu} e_{i,b}^{\mu'}\biggr).
\end{gather*}
We rewrite this as
\begin{gather*}
\sum_{\alpha\colon i \rightarrow j} \sum_{a \in \mathbb{Z}} \biggl(\sum_{k > 0} \sum_{\mu, \mu'} e_{i,a}^{\mu} e_{j,a+k+1}^{\mu'} + \sum_{\mu < \mu'} e_{i,a}^{\mu} e_{j,a+1}^{\mu'}\biggr) \\
\qquad- \sum_{i} \sum_{a \in \mathbb{Z}} \biggl(\sum_{k > 0} \sum_{\mu, \mu'} e_{i,a}^{\mu} e_{i,a+k}^{\mu'} + \sum_{\mu < \mu'} e_{i,a}^{\mu} e_{i,a}^{\mu'}\biggr) ,
\end{gather*}
by replacing $b$ with $a+k+1$ in the first term and with $a+k$ in the second. In terms of $k$-shifts, this can be rewritten as
\begin{align*}
- \sum_{k > 0} \sum_{\mu, \mu'} \big\langle\mathbf{e}_{\mu}, {\tau}^k \mathbf{e}_{\mu'}\big\rangle - \sum_{\mu < \mu'} \langle\mathbf{e}_{\mu}, \mathbf{e}_{\mu'}\rangle.
\end{align*}

As before, we consider the ring of localized motives $ \mathbb{M} = K_0 (\text{Var}_{\mathbb{C}}) \big[\mathbb{L}^{-1}, \bigl(1- \mathbb{L}^i\bigr)^{-1}, i \geq 1\big]$ and conclude from Theorem~\ref{theo mot} the following.

\begin{Proposition} \label{prop iden}
In $\mathbb{M} \llbracket Q_0 \rrbracket$, we can express the nilpotent motivic generating function
\begin{align*}
\sum_{\mathbf{d}} \frac{[ N_{\mathbf{d}}]}{[G_{\mathbf{d}}]} t^{\mathbf{d}}
\end{align*}
as
\begin{align*}
 \sum
 \mathbb{L}^{ - \sum_{k > 0} \sum_{\mu, \mu' \in [0,1] \cap \mathbb{Q}} \langle\mathbf{e}_{\mu}, {\tau}^k \mathbf{e}_{\mu'}\rangle - \sum_{\mu < \mu'} \langle\mathbf{e}_{\mu}, \mathbf{e}_{\mu'}\rangle} \prod_{\mu \in [0,1] \cap \mathbb{Q}} \frac{\big[ R_{\mathbf{e}_\mu}^{\theta\text{\rm -sst}}(Q_{\mathbb{Z}}) \big]}{[ G_{\mathbf{e}_\mu}]} t^{|\sum_{\mu} \mathbf{e}_{\mu}|} ,
\end{align*}
where the sum ranges over $L({\bf d})$, that is, over all tuples $({\bf e}_\mu)_{\mu\in[0,1]\cap\mathbb{Q}}$ of dimension vectors for~$Q_\mathbb{Z}$ such that ${\bf e}_\mu=0$ for all but finitely many $\mu$, and $\sigma({\bf e}_\mu)=\mu$ otherwise.
\end{Proposition}

\begin{Example} As a quick plausibility check, we consider the quiver $Q$ of extended Dynkin type~$\tilde{A}_2$ with three vertices $i$, $j$, $k$ and arrows $i\rightarrow j\rightarrow k$ as well as $i\rightarrow k$, with dimension vector~${{\bf d}={\bf i}+{\bf j}+{\bf k}}$ (that is, $d_i=d_j=d_k=1$). The part of $Q_\mathbb{Z}$ on which all relevant ${\bf e}_\mu$ are supported~is
\[(j,-1)\rightarrow(k,0)\leftarrow(i,-1)\rightarrow(j,0)\rightarrow(k,1)\leftarrow(i,0).\]
We have seven tuples $({\bf e}_\mu)_\mu$ such that $|\sum_\mu{\bf e}_\mu|={\bf d}$ for which semistable representations exist
\begin{gather*}
 ({\bf(i,-1)}+{\bf(j,0)}+{\bf(k,1)})_0,\qquad
({\bf(i,-1)}+{\bf(j,0)}+{\bf(k,0)})_{1/3},\\
({\bf(i,-1)}+{\bf(j,-1)}+{\bf(k,0)})_{2/3},\qquad
(({\bf(i,0)})_0,({\bf(j,-1)}+{\bf(k,0)})_{1/2}),\\
(({\bf(j,0)})_0,({\bf(i,-1)}+{\bf(k,0)})_{1/2}),\qquad
(({\bf(k,0)})_0,({\bf(i,-1)}+{\bf(j,0)})_{1/2}),\\
({\bf(i,0)}+{\bf(j,0)}+{\bf(k,0)})_0.\end{gather*}
The only non-zero $\mathbb{L}$-exponent in the above summation equals $1$ for the first tuple. Semistability in each case is given by non-vanishing of all scalars representing the arrows. This yields
\[[N_{\bf d}]=(\mathbb{L}+2)(\mathbb{L}-1)^2+3(\mathbb{L}-1)+1=\mathbb{L}^3,\]
as expected.
\end{Example}

\begin{Definition}
We define a new twisted formal power series ring $\mathbb{M}_{\mathbb{L},*}\llbracket(Q_\mathbb{Z})_0\rrbracket$ for $Q_{\mathbb{Z}}$ with topological $\mathbb{M}$-basis elements $u^{\bf e}$ for ${\bf e}$ a dimension vector for $Q_\mathbb{Z}$ with finite support, and multiplication
given by
\begin{align*}
u^{\mathbf{e}} \ast u^{\mathbf{e}'} = \mathbb{L}^{-\sum_{k\geq 0} \langle\mathbf{e}, {\tau}^k \mathbf{e}'\rangle_{\mathbb{Q}_{\mathbb{Z}}} - \sum_{k > 0} \langle\mathbf{e}', {\tau}^k \mathbf{e}\rangle_{\mathbb{Q}_{\mathbb{Z}}}} \cdot u^{\mathbf{e}+\mathbf{e}'} .
\end{align*}

We have a specialization map of $\mathbb{M}$-modules $\mathbb{M}_{\mathbb{L},*}\llbracket(Q_\mathbb{Z})_0\rrbracket\rightarrow \mathbb{M}_\mathbb{L}\llbracket Q_0\rrbracket$ by mapping $u^{\bf e}$ to $t^{|{\bf e}|}$ (in general not compatible with the multiplications).
\end{Definition}
This definition now allows us to interpret Proposition \ref{prop iden} as a product formula using an ascending ordered product over slopes. Namely, given series $F_\mu\in\mathbb{M}_{\mathbb{L},*}\llbracket(Q_\mathbb{Z})_0\rrbracket$ without constant term for $\mu\in[0,1]\cap\mathbb {Q}$ as below, we define
\[\prod_{\mu\in[0,1]\cap\mathbb{Q}}^\rightarrow(1+F_\mu)=\sum_{s\geq 0} \,\, \sum_{0\leq\mu_1<\dots<\mu_s<1}\,\, \prod_{k=1}^sF_{\mu_k}.\]
Then we find by a direct calculation the following.

\begin{Theorem} \label{theoremwcf}
The nilpotent motivic generating function of Proposition {\rm\ref{prop iden}} equals the specialization of the following ordered product in $\mathbb{M}_{\mathbb{L},*}\llbracket(Q_\mathbb{Z})_0\rrbracket$:
\begin{align*}
\prod^{\rightarrow}_{\mu \in [0,1] \cap \mathbb{Q}} \Biggl(1+ \sum\limits_{\substack{\sigma({\bf e})=\mu \\
\supp(\mathbf{e}) \ \mathrm{finite}}} \mathbb{L}^{- \sum_{k>0} \langle\mathbf{e}, {\tau}^k \mathbf{e}\rangle} \frac{\big[R_{\mathbf{e}}^{\theta\text{\rm -sst}} (Q_{\mathbb{Z}})\big]}{[G_{\mathbf{e}}]} \cdot u^{\mathbf{e}} \Biggr) .
\end{align*}
\end{Theorem}

\subsection*{Acknowledgments} The authors would like to thank anonymous referees for suggestions on improving the presentation and for pointing out \cite[Lemma 2.6]{BSV}. The second named author would like to thank S.~Mozgovoy for stimulating discussions. The work of both authors is supported by the DFG GRK 2240 ``Algebro-Geometric Methods in Algebra, Arithmetic and Topology''.

\pdfbookmark[1]{References}{ref}
\LastPageEnding

\end{document}